\documentclass[a4paper,11pt]{article}

\usepackage[utf8]{inputenc} 
\usepackage[T2A]{fontenc} 
\usepackage[russian,english]{babel}

\usepackage[strict]{chngpage}
\usepackage{adjustbox}

\usepackage{latexsym}
\usepackage{graphicx}
\usepackage{amsfonts}
\usepackage{amsthm}
\usepackage{amsmath}
\usepackage{amssymb}

\usepackage{accents}
\usepackage{indentfirst}
\usepackage{float}
\usepackage{mathtools}
\usepackage[noadjust]{cite}
\usepackage{hyperref}

\usepackage{enumitem}
\setlist[enumerate]{font=\upshape,label=(\arabic*)}
\setlist{noitemsep} 

\allowdisplaybreaks

\newtheorem{theorem}{Theorem}
\newtheorem{proposition}[theorem]{Proposition}
\newtheorem{corollary}[theorem]{Corollary}
\newtheorem{lemma}[theorem]{Lemma}
\theoremstyle{definition}
\newtheorem{definition}[theorem]{Definition}
\newtheorem{notation}[theorem]{Notation}
\theoremstyle{definition}
\newtheorem{remark}[theorem]{Remark}

\numberwithin{theorem}{section}
\numberwithin{equation}{section} 
\numberwithin{figure}{section} 
\numberwithin{table}{section} 

\usepackage{geometry} 

\newcommand{\RNumb}[1]{\uppercase\expandafter{\romannumeral #1\relax}}

\DeclareMathOperator{\lcm}{lcm}

\providecommand{\keywords}[1]{\textbf{Keywords:} #1}
\providecommand{\msc}[1]{\textbf{MSC 2020:} #1}

\newcommand\freefootnote[1]{%
    \bgroup
    \renewcommand\thefootnote{\fnsymbol{footnote}}%
    \renewcommand\thempfootnote{\fnsymbol{mpfootnote}}%
    \footnotetext[0]{#1}%
    \egroup
}

\begin{document}
\title{Linear transformations of tropical matrices\\ preserving individual values of the cyclicity index}

\author{
Alexander Guterman$^{a}$, Alexander Vlasov$^{b,c}$, and Svetlana Zhilina$^{b,c}$
}
\date{\small \em
$^a$Bar-Ilan University, Ramat Gan, 5290002, Israel\\
$^b$Lomonosov Moscow State University, Moscow, 119991, Russia\\
$^c$Moscow Center of Fundamental and Applied Mathematics, Moscow, 119991, Russia
}

 \maketitle

\begin{abstract}
The cyclicity index is used to study the periodic behavior of powers of tropical matrices. It is determined by the structure of the critical subgraph of a matrix, which consists of the cycles with the largest average weight. In this paper bijective linear transformations of tropical matrices that preserve an individual value $k\neq 2$ of this invariant are studied.
\end{abstract}


\msc{15A80, 05C38, 05C50.}

\freefootnote{The authors' research was supported by the Ministry of Science and Higher Education of the Russian Federation as part of the program of the Moscow Center of Fundamental and Applied Mathematics under agreement No. 075-15-2025-345.}

\freefootnote{Email addresses: \texttt{alexander.guterman@biu.ac.il} (Alexander Guterman), \texttt{aleksandr.vlasov@math.msu.ru} (Alexander Vlasov), \texttt{s.a.zhilina@gmail.com} (Svetlana Zhilina)}

\section{Introduction}

A set $\mathcal{S}$ with two binary operations, addition and multiplication, is called a \textit{semiring} if
\begin{enumerate}
    \item $\mathcal{S}$ is a commutative monoid under addition with identity element $0$;
    \item $\mathcal{S}$ is a semigroup under multiplication;
    \item multiplication is distributive over addition;
    \item $a\cdot 0 = 0\cdot a = 0$ for any $a\in \mathcal{S}$.
\end{enumerate}
A semiring is called \textit{commutative} if its multiplication is commutative, \textit{anti-negative} if the equality $a+b = 0$ implies $a = b = 0$, and \textit{idempotent} if $a+a = a$ for any $a\in \mathcal{S}$. Every idempotent semiring is anti-negative. Indeed, if $a+b = 0$, then $0 = a+b = a+a+b = a+0 = a$, and similarly $b = 0$.

In this paper we consider the \textit{tropical semiring}
$\mathbb{R}_{\max} = \mathbb{R}\cup \{-\infty\}$ with the operations of tropical addition $\oplus$ and tropical multiplication $\otimes$ defined as follows:
\begin{samepage}
    \begin{align*}
    &a\oplus b = \max\{a, b\},\\
    &a\otimes b = a+b
    & &\text{for all }a, b\in \mathbb{R},\\
    &a\oplus (-\infty) = (-\infty)\oplus a = a,\\
    &a\otimes (-\infty) = (-\infty)\otimes a = -\infty
    & &\text{for any } a\in \mathbb{R}_{\max}.
\end{align*}
The element $-\infty$ is the identity element with respect to tropical addition, and $0$ is the identity element with respect to tropical multiplication. This semiring is commutative and idempotent and, therefore, anti-negative.
\end{samepage}

If $\mathcal{S}$ is a semiring, then the set $M_n(\mathcal{S})$ of $n\times n$ matrices over $\mathcal{S}$ is also a semiring. For brevity, we denote the set of matrices over the tropical semiring by~$M_n$. Throughout the paper we assume that $n\geq 2$.

When studying tropical matrices, it is convenient to consider weighted directed graphs. Namely, to each matrix $A\in M_n$ we associate the graph $G(A)$ on the vertices $1, \ldots, n$ such that the arc from vertex $i$ to vertex $j$ has weight $A_{ij}$ if $A_{ij}\neq -\infty$, and is absent otherwise. All graphs under consideration have no multiple arcs but may contain loops. A matrix $A$ is called \textit{irreducible} if the graph $G(A)$ is strongly connected. By a \textit{cycle} we always mean a simple directed cycle, that is, a closed path with no repeated vertices except that the initial and the terminal vertices coincide. A loop is a cycle of length~$1$.

The main object of study in this paper is an invariant called the \textit{cyclicity index}. It is defined first for directed graphs and then for tropical matrices.

\begin{definition}
    Let $G$ be a directed graph. Its \textit{cyclicity index} $\sigma_G$ is defined as follows:
    \begin{itemize}
        \item if $G$ is strongly connected and contains no cycles, then $\sigma_G = 1$;
        \item if $G$ is strongly connected and contains at least one cycle, then $\sigma_G$ is equal to the greatest common divisor of the lengths of all cycles in~$G$;
        \item if $G$ is not strongly connected, then $\sigma_G$ is equal to the least common multiple of the cyclicity indices of all its strongly connected components.
    \end{itemize}
\end{definition}

\begin{samepage}
\begin{definition}
    Let $G$ be a weighted directed graph, and let $C$ be a cycle of length $l$ in~$G$ consisting of arcs with weights $a_1, \ldots, a_l$. The \textit{weight} of the cycle $C$ is the tropical product of the weights of its arcs:
    \[
        w(C) = a_1\otimes \ldots\otimes a_l = a_1+\ldots+a_l.
    \]
    The \textit{average weight} $a(C)$ of the cycle $C$ is the tropical geometric mean of the weights of its arcs or, equivalently, their usual arithmetic mean:
    \[
        a(C) = \sqrt[\otimes l]{a_1\otimes\ldots\otimes a_l} = \frac{a_1+\ldots+a_l}{l} = \frac{w(C)}{l}.
    \]
    Similarly, the \textit{weight} $w(G)$ of an arbitrary weighted directed graph $G$ is the tropical product of the weights of all its arcs. For a graph without arcs we set $w(G)=0$.
\end{definition}
\end{samepage}

\begin{definition}
    Let $G$ be a weighted directed graph. A cycle $C$ in $G$ is called \textit{critical} if $a(C)$ equals the largest average weight of a cycle in~$G$. The \textit{critical subgraph} $G_c$ is the union of all critical cycles in~$G$.
\end{definition}

\begin{definition}
    Let $A\in M_n$. The \textit{cyclicity index} $\sigma(A)$ of the matrix $A$ is the cyclicity index of the critical subgraph $G_c(A)$ of the graph~$G(A)$. If $G(A)$ contains no cycles, we set $\sigma(A) = 1$.
\end{definition}

As the following theorem shows, the cyclicity index plays an important role in the study of the behavior of powers of tropical matrices.

\begin{theorem}[{\cite[Theorem~3.9]{HOW}}]
    Let $A\in M_n$ be an irreducible tropical matrix with finite eigenvalue $\lambda$ and cyclicity index $\sigma$. Then there exists a positive integer~$N$ such that for any positive integer $t\geq N$ we have
    \[
        A^{\otimes(t+\sigma)} = \lambda^{\otimes \sigma}\otimes A^{\otimes t}.
    \]
\end{theorem}

In particular, if $\lambda = 0$, then, starting from some point on, the sequence of powers of~$A$ is periodic with period~$\sigma$. Thus, the cyclicity index is used for computing large powers of tropical matrices, see~\cite{HOW}. Bounds on the length of the transient phase before the periodic behavior is achieved were obtained, in particular, in~\cite{HA,CFN,KMNS}.

Once an invariant is defined, the question naturally arises as to which transformations preserve it. One of the first to study a problem of this kind was Frobenius~\cite{F}, who characterized the bijective linear mappings of complex matrices that preserve the determinant. Later this result was generalized by Schur~\cite{Schur}, who described the mappings preserving all minors of a fixed order. Dieudonn\'e~\cite{Di} obtained a characterization of the bijective linear mappings preserving the set of singular matrices over an arbitrary field.

Subsequently, many authors studied mappings that preserve various invariants of matrices over various algebraic structures, see~\cite{BP1,BP2,GM1,Pierce} and the references therein.

We say that a mapping $T: M_n\to M_n$ \textit{preserves the cyclicity index} if $\sigma(T(A)) = \sigma(A)$ for any matrix $A\in M_n$. A mapping $T$ \textit{preserves the value $k$ of the cyclicity index} if for any matrix $A\in M_n$ the equality $\sigma(A) = k$ implies $\sigma(T(A)) = k$. A mapping $T$ is called \textit{linear} if
\[
    T(\alpha\otimes A\oplus \beta\otimes B)
    = \alpha\otimes T(A)\oplus \beta\otimes T(B)
\]
for any $A, B\in M_n$ and $\alpha, \beta\in \mathbb{R}_{\max}$.

Linear mappings of tropical matrices that preserve the cyclicity index were previously considered in~\cite{GKT,GKV1,GKV2}. In~\cite{GKT} the bijective linear mappings preserving the cyclicity index were described, and in~\cite{GKV1} a complete characterization of such mappings without the bijectivity assumption was given. In~\cite{GKV2} it was shown that preserving all values of the cyclicity index is equivalent to preserving its values $1$ and~$2$ simultaneously. In addition, in~\cite{GKV3} linear mappings preserving the cyclicity index were studied for matrices over the semirings of nonnegative real numbers and nonnegative integers with the usual addition and multiplication.

In this paper we consider bijective linear mappings of tropical matrices that preserve a single fixed value $k \neq 2$ of the cyclicity index. To describe such mappings, we use the results on additive transformations of unweighted directed graphs preserving individual values of cycle lengths, obtained in our previous paper~\cite{GVZ1}.

The paper is organized as follows. Section~2 contains some definitions and known facts on bijective linear transformations of $M_n$ and their connections with additive mappings of unweighted directed graphs. In Section~3 we obtain a complete characterization of the bijective linear mappings of tropical matrices that preserve the value $1$ of the cyclicity index, see Theorem~\ref{main_theorem}. In Corollary~\ref{cor_preserving_cyclicity} we show that if a mapping and its inverse both preserve this value, then they preserve all values of the cyclicity index. Section~4 is devoted to bijective linear mappings preserving a fixed value \(k\geq3\) of the cyclicity index. In Theorem~\ref{th_non_diagonal_elements_k} we describe their action on the off-diagonal matrix units, and in Theorem~\ref{th_final_k} we obtain a complete characterization of such mappings for $n\geq\max(k,4)$. Corollary~\ref{cor_preserving_k} shows that if $k$ and $n$ satisfy certain restrictions, then preserving the value~$k$ of the cyclicity index is equivalent to preserving all values, and Corollary~\ref{cor_preserving_k_both_sides} is an analogue of Corollary~\ref{cor_preserving_cyclicity} for $k \geq 3$.

\section{Definitions and known results}

For convenience, we introduce the following notation. The \textit{matrix unit} $E_{ij}\in M_n$ is the matrix whose $(i, j)$-th entry equals $0$, while all other entries are equal to~$-\infty$. A matrix $P\in M_n$ is called a \textit{permutation matrix} if there exists a permutation $\tau\in S_n$ such that $P_{i, \tau(i)} = 0$ for all $i = 1, \ldots, n$, while all other entries are equal to $-\infty$. The \textit{Schur product} of matrices $A, B\in M_n$ is the matrix $A\circ B$ defined by the equalities
\[
    (A\circ B)_{ij} = A_{ij}\otimes B_{ij}, \qquad 1\leq i, j\leq n.
\]
We denote by $\mathcal{O}_n$ the $n\times n$ matrix all of whose entries are equal to~$0$.

A special role is played by matrices which correspond to graphs with a single arc. We denote
    \begin{align*}
        \mathcal{E}_n &= \{\alpha\otimes E_{ij}\mid 1\leq i, j\leq n, \ \alpha\in \mathbb{R}\}, \\
        \mathcal{D}_n &= \{\alpha\otimes E_{ii}\mid 1\leq i\leq n, \ \alpha\in \mathbb{R}\}, \\
        \mathcal{N}_n &= \{\alpha\otimes E_{ij}\mid 1\leq i \neq j\leq n, \ \alpha\in \mathbb{R}\}.
    \end{align*}
Thus, $\mathcal{E}_n = \mathcal{D}_n\cup \mathcal{N}_n$, where $\mathcal{D}_n$ consists of the matrices whose graph has a single arc which is a loop, and $\mathcal{N}_n$ consists of the matrices for which this arc is not a loop.

Let $\widehat M_n$ be the set of unweighted directed graphs without multiple arcs on the fixed vertex set $\{1, \ldots, n\}$. For unweighted graphs we use notation similar to that which was introduced for tropical matrices. We denote by $\mathcal{\widehat E}_n$ the set of graphs from $\widehat M_n$ with a single arc, by $\mathcal{\widehat D}_n$ its subset consisting of the graphs with a single loop, and by $\mathcal{\widehat N}_n$ the subset of graphs with a single non-loop arc. We identify the graph with a single arc from vertex $i$ to vertex $j$ with this arc and denote it by $e_{ij}$.

For graphs $\widehat G_1, \widehat G_2$ from $\widehat M_n$ we define their union $\widehat G_1\cup \widehat G_2$ as the union of the arc sets of these graphs. We call a mapping $S: \widehat M_n\to \widehat M_n$ \textit{additive} if
\[
    S(\widehat G_1\cup \widehat G_2)
    = S(\widehat G_1)\cup S(\widehat G_2)
\]
for any $\widehat G_1, \widehat G_2\in \widehat M_n$.

Apart from the weighted graph $G(A)$, to each matrix $A\in M_n$ we can also associate the unweighted directed graph $\widehat G(A)$ obtained from $G(A)$ by deleting the weights of all arcs. Note that for any matrices $A, B \in M_n$ and $\alpha \in \mathbb{R}$ we have $\widehat G(A \oplus B) = \widehat G(A) \cup \widehat G(B)$ and $\widehat G(\alpha \otimes A) = \widehat G(A)$. Therefore, every linear mapping of matrices $T: M_n\to M_n$ induces an additive mapping of graphs $\widehat{T}: \widehat M_n \to \widehat M_n$ given by the formula $\widehat{T}(\widehat G(A)) = \widehat G(T(A))$ for all $A \in M_n$. Namely, we define $\widehat T(e_{ij})$ as the graph $\widehat{G}(T(E_{ij}))$ and extend the mapping~$\widehat T$ by additivity to the whole set $\widehat M_n$. From the following lemma it is easy to deduce that if $T$ is bijective, then $\widehat{T}$ is bijective as well, since in this case $\widehat T$ is a bijection on $\mathcal{\widehat E}_n$.

\begin{lemma}[{\cite[Lemma~3.1]{GKT}}]\label{surjective_maps}
    Let $T: M_n\to M_n$ be a linear mapping of tropical matrices. Then the following conditions are equivalent:
    \begin{enumerate}
        \item $T$ is bijective;
        \item $T$ is surjective;
        \item $T$ is a bijection on the set $\mathcal{E}_n$.
    \end{enumerate}
\end{lemma}

\begin{definition}\label{def_cycle_length_preserving}
    We say that a mapping $S: \widehat M_n\to \widehat M_n$ \textit{preserves the value $k$ of cycle lengths} if it maps every graph whose arcs are exactly the arcs of a single cycle of length~$k$ to a graph of the same form. A mapping $S$ \textit{preserves cycle lengths} if it preserves all their values.
\end{definition}

In Lemmas~\ref{lem_cycle_to_cycle} and~\ref{lem_preserving_k_cycles} below we show that:
\begin{enumerate}
    \item if a bijective linear mapping $T: M_n\to M_n$ preserves the value $1$ of the cyclicity index, then the induced mapping $\widehat T:~\widehat M_n\to\widehat M_n$ preserves cycle lengths;
    \item if a bijective linear mapping $T: M_n\to M_n$ preserves the value $k$ of the cyclicity index, $k\leq n$, then the induced mapping $\widehat T$ preserves the value $k$ of cycle lengths.
\end{enumerate}

Thus, to describe the bijective linear mappings of tropical matrices preserving individual values of the cyclicity index, we need the characterization of additive mappings of unweighted directed graphs preserving individual values of cycle lengths, which was obtained in~\cite{GVZ1}. We recall the main results of that paper.

\begin{proposition}[{\cite[Proposition~3.1]{GVZ1}}]\label{prop_bij_action_on_k_cycles}
    Let $S:\widehat M_n\to \widehat M_n$ be an additive mapping of unweighted graphs such that the restriction $S|_{\mathcal{\widehat N}_n}$ is bijective, and let $2\leq k\leq n$. The mapping $S$ preserves the value~$k$ of cycle lengths if and only if it maps the set of cycles of length~$k$ bijectively onto itself.
\end{proposition}

\begin{lemma}[{\cite[Lemma~2.3]{GVZ1}}]\label{lem_vertices_permutation}
Let $S:\widehat M_n\to \widehat M_n$ be an additive mapping of unweighted graphs preserving the values $2$ and~$3$ of cycle lengths, $n \geq 2$, and let the set $\mathcal{\widehat N}_n$ be invariant under the mapping~$S$. Then
    \begin{enumerate}
        \item the restriction $S|_{\mathcal{\widehat N}_n}$ is bijective;
        \item there exists a permutation $\tau\in S_n$ such that exactly one of the following two conditions holds simultaneously for all $i\ne j$:
            \[
            S(e_{ij})=e_{\tau(i)\tau(j)} \quad \text{or} \quad S(e_{ij})=e_{\tau(j)\tau(i)}.
            \]
        In other words, if a graph $\widehat G$ contains no loops, then $S(\widehat G)$ is obtained from $\widehat G$ by a renumbering of the vertices and, possibly, by reversing the orientation of all arcs.
        \end{enumerate}
\end{lemma}

\begin{theorem}[{\cite[Theorem~4.5]{GVZ1}}]\label{th_vertices_permutation}
    Let $S:\widehat M_n\to \widehat M_n$ be an additive mapping of unweighted graphs preserving the value $k\geq 3$ of cycle lengths, $n\geq\max(k,4)$, and let the set $\mathcal{\widehat N}_n$ be invariant under the mapping $S$. 
    \begin{enumerate}
        \item If $n \geq \max(k+1,5)$, then the restriction $S|_{\mathcal{\widehat N}_n}$ is bijective.
        \item If the restriction $S|_{\mathcal{\widehat N}_n}$ is bijective, then there exists a permutation $\tau\in S_n$ such that exactly one of the following two conditions holds simultaneously for all $i\ne j$:
        \[
        S(e_{ij})=e_{\tau(i)\tau(j)} \quad \text{or} \quad S(e_{ij})=e_{\tau(j)\tau(i)}.
        \]
        \end{enumerate}
\end{theorem}

\section{Bijective linear mappings of tropical matrices preserving the value $1$ of the cyclicity index} 

By Lemma~\ref{surjective_maps}, if $T$ is a bijective linear mapping of tropical matrices, then we may assume that for all $1 \leq i,j \leq n$ it holds that
\begin{equation} \label{equation:basis-images}
    T(E_{ij}) = \alpha_{ij}\otimes E_{km}
\end{equation}
for some $\alpha_{ij} \in \mathbb{R}$ and $1\leq k,m\leq n$.

\begin{proposition}\label{prop_diagonal_to_diagonal}
    Let $T: M_n\to M_n$ be a bijective linear mapping of tropical matrices of order $n$ preserving the value $1$ of the cyclicity index. Then for any matrix $A\in\mathcal{D}_n$ we have $T(A)\in\mathcal{D}_n$.
\end{proposition}

\begin{proof}
    Assume the contrary. Let $A = \alpha\otimes E_{ll}$ and $T(A)\notin\mathcal{D}_n$. By Lemma~\ref{surjective_maps}, we have $T(A)\in\mathcal{N}_n$. Then $T(A)=\alpha\otimes\alpha_{ll}\otimes E_{km}$, where $k\neq m$. 

    Let $\lambda$ be the largest average weight of cycles in $T(\mathcal{O}_n)$, and let $\mu$ be the average weight of the cycle $k\to m\to k$ in $T(\mathcal{O}_n)$.
    Consider the matrix 
    \[
    B=\mathcal{O}_n\oplus (6(\lambda-\mu)+1)\otimes E_{ll}.
    \]
    Since $\lambda\geq\mu$ by definition, the weight of the loop at the vertex $l$ in the graph $G(B)$ is at least $1$. Moreover, the average weights of all other cycles in $G(B)$ are equal to~$0$; hence $\sigma(B) = 1$.

    By linearity, we have
    \begin{align*}
    T(B) &= T(\mathcal{O}_n)\oplus (6(\lambda-\mu)+1)\otimes T(E_{ll}) \\
    &= T(\mathcal{O}_n)\oplus(6(\lambda-\mu)+1)\otimes\alpha_{ll}\otimes E_{km}.
    \end{align*}
    Since $(T(\mathcal{O}_n))_{km} = \alpha_{ll}$ and $6(\lambda-\mu)+1 > 0$, when passing from $G(T(\mathcal{O}_n))$ to $G(T(B))$ the weight of the arc $k\to m$ increases by $6(\lambda-\mu)+1$, while the weights of all other arcs do not change.
    Then the weights of the loops in $G(T(B))$ do not exceed~$\lambda$, the average weight of the cycle $k\to m\to k$ is equal to
    \[
    \mu + \tfrac{6(\lambda-\mu)+1}{2} = 3\lambda - 2\mu + \tfrac{1}{2}\geq \lambda + \tfrac{1}{2},
    \]
    and the average weight of any other cycle does not exceed
    \[
    \lambda+\tfrac{6(\lambda-\mu)+1}{3} = 3\lambda-2\mu+\tfrac{1}{3}.
    \]
    Thus the critical subgraph $G_c(T(B))$ consists of the single cycle $k\to m\to k$, and $\sigma(T(B)) = 2$, which contradicts the fact that the mapping $T$ preserves the value $1$ of the cyclicity index.
\end{proof}

\begin{corollary}\label{bijective_on_diagonal_cells_cor}
    Let $T: M_n\to M_n$ be a bijective linear mapping of tropical matrices preserving the value $1$ of the cyclicity index. Then the mappings $T|_{\mathcal{N}_n}, T|_{\mathcal{D}_n}, \widehat T|_{\mathcal{\widehat N}_n}, \widehat T|_{\mathcal{\widehat D}_n}$ are well defined and bijective.
\end{corollary}

\begin{proof}
    Proposition~\ref{prop_diagonal_to_diagonal} implies that the mapping $T|_{\mathcal{D}_n}$ is well defined. By Lemma~\ref{surjective_maps}, the mapping $T$ acts bijectively on the set $\mathcal{E}_n$. Therefore, $T|_{\mathcal{N}_n}$ is also well defined, and these restrictions are bijective. This directly implies that the mappings $\widehat T|_{\mathcal{\widehat N}_n}, \widehat T|_{\mathcal{\widehat D}_n}$ are well defined and bijective.
\end{proof}

\begin{lemma}\label{lem_cycle_to_cycle}
    Let $T: M_n\to M_n$ be a bijective linear mapping of tropical matrices preserving the value $1$ of the cyclicity index. Then the mapping $\widehat T:~\widehat M_n\to\widehat M_n$ preserves cycle lengths.
\end{lemma}
\begin{proof}
    Assume the contrary. Consider the smallest $k \in \mathbb{N}$ such that the mapping $\widehat T$ does not preserve the value $k$ of cycle lengths. By Corollary~\ref{bijective_on_diagonal_cells_cor}, $\widehat T$ acts bijectively on loops and on non-loop arcs. Hence $k \geq 2$, and the mapping~$\widehat T^{-1}$ is well defined. We show that $\widehat T^{-1}$ does not preserve the value~$k$ of cycle lengths. Indeed, otherwise, by Proposition~\ref{prop_bij_action_on_k_cycles}, it would map the set of cycles of length $k$ bijectively onto itself, but then $\widehat T$ would also preserve the value $k$ of cycle lengths. Hence there exists a cycle $\widehat C_k\in \widehat M_n$ of length~$k$ such that $\widehat T^{-1}(\widehat C_k)$ is not a cycle of length $k$.

    Note that the graph $\widehat T^{-1}(\widehat C_k)$ must contain at least one cycle. Indeed, let $A\in M_n$ be a matrix such that $\widehat G(A) = \widehat T^{-1}(\widehat C_k)$. Then $\widehat G(T(A)) = \widehat C_k$. If $\widehat G(A)$ has no cycles, then
    \[
    \sigma(A) = 1 < 2\leq k = \sigma(T(A)),
    \]
    which contradicts the fact that the mapping $T$ preserves the value $1$ of the cyclicity index.

    Note that $\widehat T^{-1}(\widehat C_k)$ cannot contain cycles of length less than~$k$, since for any $1 \leq l < k$ the mapping $\widehat T$ preserves the value~$l$ of cycle lengths.
    
    Moreover, by Corollary~\ref{bijective_on_diagonal_cells_cor}, the mapping $\widehat{T}$ is a bijection on $\mathcal{\widehat N}_n$. Hence the graph $\widehat T^{-1}(\widehat C_k)$ has $k$ arcs. Therefore $\widehat T^{-1}(\widehat C_k)$ contains no cycles of length greater than $k$, and the only remaining possibility is that $\widehat T^{-1}(\widehat C_k)$ contains a cycle of length~$k$. But these two graphs have the same number of arcs, so they must coincide, a contradiction.
\end{proof}

\begin{corollary}\label{cor_adj_edges}
    Let $T: M_n\rightarrow M_n$ be a bijective linear transformation of tropical matrices preserving the value $1$ of the cyclicity index. Then there exists a permutation $\tau\in S_n$ such that exactly one of the following two conditions holds simultaneously for all $1\le i\ne j\le n$:
    \[
    T(E_{ij})=\alpha_{ij}\otimes E_{\tau(i)\tau(j)} \quad \text{or} \quad T(E_{ij})=\alpha_{ij}\otimes E_{\tau(j)\tau(i)}, \quad \alpha_{ij} \in \mathbb{R}.
    \]
\end{corollary}

\begin{proof}
    The statement follows immediately from Lemmas~\ref{lem_vertices_permutation} and~\ref{lem_cycle_to_cycle}.
\end{proof}

In~\cite{GKV2}, Lemma~2.12 was stated for linear transformations of tropical matrices preserving the values~$1$ and~$2$ of the cyclicity index. However, its proof is valid verbatim for bijective linear transformations preserving the value $1$ of the cyclicity index, see Lemma~\ref{lem_the_same_add_to_average_weight}. Indeed, by Corollary~\ref{cor_adj_edges}, every such transformation acts on non-loop arcs by a renumbering of the vertices and, possibly, reversing the orientation of all arcs. Hence the lengths of cycles, their number, and their common arcs are preserved. The rest of the argument in the proof of \cite[Lemma~2.12]{GKV2} uses only linearity and the preservation of the value~$1$ of the cyclicity index.

\begin{lemma}[{\cite[Lemma~2.12]{GKV2}}]\label{lem_the_same_add_to_average_weight}
    Let $T: M_n\rightarrow M_n$ be a bijective linear transformation of tropical matrices preserving the value $1$ of the cyclicity index. Then there exists a number $\lambda\in\mathbb{R}$ with the following property: for all $C\in M_n$ such that $G(C)$ is a cycle of length at least~$2$ it holds that
    \[
    a(G(T(C))) = \lambda\otimes a(G(C)),
    \]
    where $a(G(C))$ is the average weight of the cycle $G(C)$.
\end{lemma}

\begin{lemma}\label{lem_diagonal_permutation}
    Let $T:M_n\to M_n$ be a bijective linear transformation of tropical matrices preserving the value $1$ of the cyclicity index, let $\tau\in S_n$ be the permutation defined in Corollary~\ref{cor_adj_edges} and $\lambda\in\mathbb{R}$ be the number defined in Lemma~\ref{lem_the_same_add_to_average_weight}. Then there exist numbers $\lambda_i\in\mathbb{R}_{\geq0}$, $1\leq i\leq n$, and a permutation $\pi:\{1,2,\ldots,n\}\to\{1,2,\ldots,n\}$ such that 
    \[
    T(E_{ii}) = \lambda\otimes\lambda_i\otimes E_{(\tau\circ\pi)(i)(\tau\circ\pi)(i)}
    \]
    for all $1\leq i\leq n$. Moreover, if $n\geq 3$ and $\lambda_i = 0$, then $\pi(i) = i$.
\end{lemma}

\begin{proof}
    Proposition~\ref{prop_diagonal_to_diagonal} implies that $T(E_{ii})\in\mathcal{D}_n$. Define the permutation $\pi':\{1,\ldots,n\}\to\{1,\ldots,n\}$ so that $T(E_{ii}) = \alpha_{ii}\otimes E_{\pi'(i)\pi'(i)}$.

    We first show that $\alpha_{ii}\geq\lambda$ for all $1\leq i\leq n$. Consider the matrix 
    \[A = \bigoplus\limits_{m=1}^{n-1} E_{m,m+1}\oplus E_{n1}\oplus E_{ii}.\]
    The graph $G(A)$ consists of the cycle $1\to 2\to\ldots\to n\to 1$ and the loop at the vertex~$i$, and all arcs in $G(A)$ have weight $0$. Hence $\sigma(A)=1$. By Lemma~\ref{lem_cycle_to_cycle}, $G(T(A))$ also consists of a cycle of length $n$ and one loop of weight $\alpha_{ii}$. By Lemma~\ref{lem_the_same_add_to_average_weight}, the cycle of length $n$ has average weight $\lambda$. Since the mapping $T$ preserves the value $1$ of the cyclicity index, we have $\sigma(T(A))=1$. Hence $\alpha_{ii}\geq\lambda$. Thus for all $1\leq i\leq n$ we have $\alpha_{ii} = \lambda\otimes\lambda_i$, where $\lambda_i\geq0$.

    Now let $n\geq3$ and $\lambda_i = 0$ for some index $i$. We prove that $\pi'(i) = \tau(i)$. Assume the contrary. Consider the index $k=\tau^{-1}(\pi'(i))$ and a matrix $C\in M_n$ such that $G(C)$ is a cycle of length $n-1$ passing through all vertices except $k$, and all its arcs have weight~$0$. Set
    \[
    B = C\oplus E_{ii}.
    \]
    
    Then $G(B)$ contains a cycle of length $n-1$ passing through all vertices except the vertex $k$, and a loop at the vertex $i$; the weight of every arc is equal to $0$. Since $\pi'(i)\neq\tau(i)$, we have $k\neq i$. Hence $\sigma(B) = 1$. 
    
    At the same time, $G(T(B))$ contains a cycle of length $n-1$ passing through all vertices except the vertex $\tau(k)$, and a loop at the vertex $\pi'(i)$. By the definition of $k$, we have $\tau(k)=\pi'(i)$, therefore these two cycles do not intersect. The average weight of the cycle of length $n-1$ in the graph $G(T(B))$ is equal to~$\lambda$ by Lemma~\ref{lem_the_same_add_to_average_weight}. Since $\alpha_{ii} = \lambda\otimes\lambda_i$ and $\lambda_i=0$, the weight of the loop at the vertex $\pi'(i)$ is also equal to $\lambda$. Then the critical subgraph $G_c(T(B))$ contains two disjoint cycles of lengths $n-1$ and $1$, so $\sigma(T(B)) = n-1\neq 1$. This contradicts the fact that the mapping~$T$ preserves the value $1$ of the cyclicity index.

    Finally, set $\pi = \tau^{-1}\circ\pi'$. Then
    \[
    T(E_{ii}) = \lambda\otimes\lambda_i\otimes E_{\pi'(i)\pi'(i)} = \lambda\otimes\lambda_i\otimes E_{(\tau\circ\pi)(i)(\tau\circ\pi)(i)}.
    \]
    If $n\geq3$ and $\lambda_i = 0$, then $\pi(i) = (\tau^{-1}\circ\pi')(i) = \tau^{-1}(\pi'(i)) = \tau^{-1}(\tau(i)) = i$. 
\end{proof}

\begin{theorem}\label{main_theorem}
    Let $T:M_n\to M_n$ be a bijective linear transformation of tropical matrices. 
    Then the following conditions are equivalent:
    \begin{enumerate}
        \item the mapping $T$ preserves the value $1$ of the cyclicity index,
        \item 
        there exist
        \begin{itemize}
            \item a matrix $B = (b_{ij})\in M_n$ such that $b_{ij}\neq-\infty$ for all $1\leq i,j\leq n$, and all cycles in $G(B)$ have weight~$0$,
            \item a number $\lambda\in\mathbb{R}$,
            \item a diagonal matrix $\Lambda = \mathrm{diag}(\lambda_1,\ldots, \lambda_n)$, where all numbers $\lambda_1,\ldots,\lambda_n\in\mathbb{R}_{\geq0}$,
            \item two permutation matrices $P_1, P_2\in M_n$ with the property that if $n\geq3$ and $\lambda_i=0$ for some $1\leq i\leq n$, then $(P_2)_{ii} = 0$,
        \end{itemize}
        such that 
        \[
        T = T_1\circ T_2\circ T_3\circ T_4\circ T_5\quad\text{or}\quad T = T_2\circ T_3\circ T_4\circ T_5
        \]
        for the mappings $T_j$ defined by the formulas
        \begin{gather*}
            \begin{aligned}
            &T_1(A) = A^t, & T_2(A) = P_1\otimes A\otimes P_1^t,&\\
            &T_3(A) = \lambda\otimes A, & T_4(A) = A\circ B,&
            \end{aligned}\\
            T_5(A) = N(A)\oplus (P_2\otimes \Lambda\otimes  D(A)\otimes P_2^t)
        \end{gather*}
        for all $A\in M_n$, where $D(A)\oplus N(A) = A$, the matrix $D(A)$ is diagonal, and all diagonal entries of the matrix $N(A)$ are equal to $-\infty$.
    \end{enumerate}
\end{theorem}

\begin{proof}
    First we prove the implication $(2)\Rightarrow (1)$. To this end, we show that each of the mappings $T_1,\ldots, T_5$ is bijective and preserves the value $1$ of the cyclicity index.
    
    The bijectivity of $T_1, \ldots, T_4$ is obvious. We show that the mapping~$T_5$ is also bijective. Note that the mapping $T'_5:M_n\to M_n$ given by the equality 
    \[
    T'_5(A) =  P_2\otimes \Lambda\otimes  A\otimes P_2^t \quad \text{for all } A\in M_n
    \]
    is bijective, and the set of diagonal matrices is invariant under $T'_5$. It follows that $T_5$ does not change the off-diagonal entries of a matrix, that is,
    \[
        N(T_5(A)) = N(A), \quad D(T_5(A)) = P_2\otimes \Lambda\otimes  D(A)\otimes P_2^t,
    \]
    therefore $T_5$ is also bijective.

    Further, it is clear that $T_1, T_2, T_3$ preserve the cyclicity index. The mapping $T_4$ does not change the weights of cycles and does not add new arcs. Hence it preserves the critical subgraph $\widehat G_c(A)$ up to the weights of arcs. Thus the transformation $T_4$ also preserves the cyclicity index. 

    We show that the mapping $T_5$ preserves the value $1$ of the cyclicity index. First consider the case $n=2$. Let $X$ be a matrix such that $\sigma(X)=1$. If the graph $G(X)$ does not contain a cycle of length~$2$, then the graph $G(T_5(X))$ does not contain a cycle of length~$2$ either, since $T_5$ does not change the off-diagonal entries. Hence $\sigma(T_5(X))=1$.
    
    Now let $G(X)$ contain a cycle of length $2$. Denote its average weight by $\mu$. In this case the condition $\sigma(X)=1$ is equivalent to the inequality
    \[
    X_{11}\oplus X_{22}\geq\mu.
    \]
    
    The mapping $T_5$ permutes the loops and changes their weights by nonnegative values. Therefore
    \[
    (T_5(X))_{11}\oplus(T_5(X))_{22} \geq X_{11}\oplus X_{22}\geq\mu.
    \]
    At the same time, the average weight of the cycle of length $2$ does not change. Hence $\sigma(T_5(X))=1$.
    
    Now let $n\geq 3$. Note that the graph $G(T_5(A))$ is obtained from the graph $G(A)$ by changing the weights of loops and permuting them. Moreover, for each loop exactly one of the following two conditions holds:
    \begin{enumerate}[label=(\alph*)]
        \item the loop does not change, that is, it keeps its weight and stays in place;
        \item the weight of the loop increases, and, possibly, the loop is permuted with other loops.
    \end{enumerate}

    Consider an arbitrary matrix $X$ such that $\sigma(X)=1$. If $G(X)$ has no loops, then all diagonal entries of $X$ are equal to $-\infty$, and $T_5(X)=X$. Otherwise, $G(X)$ contains at least one cycle, and three cases are possible:
    \begin{enumerate}[label=(\roman*)]
        \item the graph $G_c(T_5(X))$ coincides with $G_c(X)$;
        \item the graph $G_c(T_5(X))$ differs from $G_c(X)$ by the appearance of new loops;
        \item the graph $G_c(T_5(X))$ consists only of loops.
    \end{enumerate}

    Indeed, the non-loop arcs in the graphs $G(X)$ and $G(T_5(X))$ coincide, while the loops may increase their weights. Then either their weights remain less than the largest average weight~$\mu$, or the weights of some of them become equal to~$\mu$, or the weights of some of them exceed~$\mu$, which corresponds to cases (i), (ii) and~(iii), respectively.

    Clearly, in cases (i) and (iii) we have $\sigma(T_5(X)) = 1$. Consider case~(ii). When a new loop is added to the critical subgraph, the partition of the vertices into strongly connected components does not change. Moreover, the cyclicity index of the component to which the loop has been added becomes equal to $1$. Hence the least common multiple of the cyclicity indices of the strongly connected components has not increased and therefore remains equal to $1$. Thus $\sigma(T_5(X))=1$, and the mapping $T_5$ preserves the value $1$ of the cyclicity index. The implication $(2)\Rightarrow (1)$ is proved.

    Now we prove the converse implication. By Corollary~\ref{bijective_on_diagonal_cells_cor}, the restrictions of $T$ to $\mathcal{N}_n$ and to $\mathcal{D}_n$ are bijective. Hence, by linearity, we may consider $T$ separately on these sets. 
    
    It follows from Corollary~\ref{cor_adj_edges} that one of the equalities $T(E_{ij})=\alpha_{ij}\otimes E_{\tau(i)\tau(j)}$ or $T(E_{ij})=\alpha_{ij}\otimes E_{\tau(j)\tau(i)}$ holds for all $1\le i\ne j\le n$ simultaneously, where $\alpha_{ij} \in \mathbb{R}$. Let the matrix $P_1$ correspond to the permutation $\tau^{-1}$. 
    Consider the mapping $T'$ such that $T = T_2\circ T'$ or $T = T_1\circ T_2\circ T'$, respectively. Then for all $1\leq i \neq j\leq n$ we obtain $T'(E_{ij}) = \alpha_{ij}\otimes E_{ij}$. In other words, the mapping $T$ acts on $\mathcal{N}_n$ as the composition of multiplication of the matrix entries by certain coefficients, $T_2$, and, possibly, $T_1$. 

    By Lemma~\ref{lem_the_same_add_to_average_weight}, the average weight of every cycle in $G(T'(A))$ is changed by some number $\lambda \in \mathbb{R}$. Define \[
    b_{ij} = 
    \begin{cases}
        \alpha_{ij}\otimes\lambda^{\otimes-1}, & i \neq j,\\
        0, & i = j.
    \end{cases}
    \]
    We show that all cycles in $G(B)$ have weight $0$. Consider an arbitrary cycle
    \[
        v_1\to v_2\to\ldots\to v_s\to v_1, \quad s\ge 2,
    \]
    where $v_i\neq v_j$ for all $i\neq j$. Let
    \[
        C=
        E_{v_1v_2}\oplus E_{v_2v_3}\oplus\ldots\oplus E_{v_sv_1}.
    \]
    Then $a(G(C))=0$, so $a(G(T'(C)))=\lambda$.
    On the other hand,
    \[
        a(G(T'(C)))=
        \frac{
            \alpha_{v_1v_2}+\alpha_{v_2v_3}+\ldots+\alpha_{v_sv_1}
        }{s}.
    \]
    Hence $\alpha_{v_1v_2}+\alpha_{v_2v_3}+\ldots+\alpha_{v_sv_1} =s\lambda$, and we obtain
    \[
        b_{v_1v_2}+b_{v_2v_3}+\ldots+b_{v_sv_1}=(\alpha_{v_1v_2}-\lambda)+(\alpha_{v_2v_3}-\lambda)+\ldots+(\alpha_{v_sv_1}-\lambda)=0.
    \]
    Thus every simple cycle of length at least $2$ in the graph $G(B)$ has weight~$0$. Loops also have weight $0$, since $b_{ii}=0$.

    Set $T_3(A) = \lambda\otimes A$ and $T_4(A) = A\circ B$. Then $T' = T_3\circ T_4\circ T''$, where $T''$ is a mapping such that $T''(E_{ij}) = E_{ij}$ for all $1\leq i\neq j\leq n$.
    Lemma~\ref{lem_diagonal_permutation} implies that $T''(E_{ii}) = \lambda_i \otimes E_{\pi(i)\pi(i)}$, and if $n\geq3$ and $\lambda_i = 0$, then $\pi(i) = i$. Let the matrix~$P_2$ correspond to the permutation $\pi^{-1}$, and let $\Lambda = \mathrm{diag}(\lambda_1,\ldots,\lambda_n)$. Then $T''(E_{ii}) = P_2\otimes\Lambda\otimes E_{ii}\otimes P_2^t$ for all $1\leq i\leq n$. Since $T''(E_{ij}) = E_{ij}$ for every pair of indices $i \neq j$, we obtain $T'' = T_5$.
\end{proof}

\begin{corollary}\label{cor_preserving_cyclicity}
    Let $T:M_n\to M_n$ be a bijective linear transformation of tropical matrices. Then the following conditions are equivalent:
    \begin{enumerate}
        \item the mappings $T$ and $T^{-1}$ preserve the value $1$ of the cyclicity index,
        \item the mappings $T$ and $T^{-1}$ preserve the cyclicity index.
    \end{enumerate}
\end{corollary}

\begin{proof}
    The implication $(2)\Rightarrow (1)$ is trivial; we prove the implication $(1)\Rightarrow (2)$.

    Since $T$ preserves the value $1$ of the cyclicity index, by Theorem~\ref{main_theorem} we have 
    \[
    T = T_1\circ T_2\circ T_3\circ T_4\circ T_5\quad\text{or}\quad T = T_2\circ T_3\circ T_4\circ T_5.
    \]
    Let $U = T_1\circ T_2\circ T_3\circ T_4$ or $U = T_2\circ T_3\circ T_4$, respectively. Then
    \[
    T = U\circ T_5.
    \]
    
    Clearly, $U$ and $U^{-1}$ preserve the cyclicity index, since the mappings $T_1,\ldots, T_4$, as well as their inverse mappings $T_1^{-1}, \ldots, T_4^{-1}$, preserve the cyclicity index.
    
    Consider $T_5^{-1}$. The equality $T = U\circ T_5$ implies that $T_5^{-1} = T^{-1}\circ U$. Since the mappings $T^{-1}$ and $U$ preserve the value 1 of the cyclicity index, $T_5^{-1}$ also preserves it.
    
    Let $\Lambda, P_2$ be the matrices from the definition of $T_5$, where $P_2$ corresponds to the permutation $\pi^{-1}$. Denote $-\Lambda = \mathrm{diag}(-\lambda_1, \dots, -\lambda_n) = \mathrm{diag}(\lambda_1^{\otimes -1}, \dots, \lambda_n^{\otimes -1})$. One can verify that
    \[
    T_5^{-1}(A) = N(A) \oplus ((-\Lambda)\otimes P_2^t\otimes D(A)\otimes P_2)
    \]
    for all $A \in M_n$, that is,
    \[
    T_5^{-1}(E_{ij}) = 
    \begin{cases}
        E_{ij}, & i \neq j,\\
        (-\lambda_{\pi^{-1}(i)})\otimes E_{\pi^{-1}(i)\pi^{-1}(i)}, & i = j.
    \end{cases}
    \]
    It then follows from Lemma~\ref{lem_diagonal_permutation} applied to $T_5^{-1}$ with the identity permutation $\tau$ and $\lambda = 0$ that $-\lambda_{\pi^{-1}(i)} \geq 0$ for all $1 \leq i \leq n$. Therefore all $\lambda_i = 0$, and $\Lambda = I_n$ is the tropical identity matrix.
    
    By the definition of $P_2$, for $n\geq 3$ we readily obtain that $T_5$ is the identity mapping. If $n=2$, then either $T_5$ is the identity mapping, or it swaps the loops $E_{11}$ and $E_{22}$, leaving the off-diagonal entries in place, that is, 
    \[
    T_5
    \begin{pmatrix}
        a& b\\
        c& d
    \end{pmatrix}
    = 
    \begin{pmatrix}
        d& b\\
        c& a
    \end{pmatrix}.
    \]
    In both cases the mappings $T_5$ and $T_5^{-1}$ preserve the cyclicity index. Hence $T$ and $T^{-1}$ preserve the cyclicity index.
\end{proof}

\section{Bijective linear mappings of tropical matrices preserving the value $k\geq3$ of the cyclicity index}

\subsection{Action on the off-diagonal entries of matrices}

\begin{lemma}\label{lem_preserving_k_cycles}
    Let $T: M_n\to M_n$ be a bijective linear mapping of tropical matrices preserving the value $k$ of the cyclicity index, $k\leq n$. Then the induced mapping $\widehat T$ preserves the value $k$ of cycle lengths.
\end{lemma}

\begin{proof}
    For $k=1$ the statement follows from Proposition~\ref{prop_diagonal_to_diagonal}.

    Now let $k>1$, and let a matrix $A\in M_n$ be such that $G(A)$ is a cycle of length $k$. By Lemma~\ref{surjective_maps}, $T(A)$ contains exactly~$k$ arcs. Moreover, since the mapping $T$ preserves the value~$k$ of the cyclicity index, we have $\sigma(T(A)) = \sigma(A) = k$. Hence $G(T(A))$ contains at least 1 cycle.

    Suppose that the graph $G(T(A))$ contains a cycle of length $d<k$. Let 
    \[
    A = \beta_1\otimes E_1\oplus \beta_2\otimes E_2\oplus\ldots\oplus\beta_k\otimes E_k, \qquad \beta_i\in\mathbb{R}, \ E_i\in\mathcal{N}_n, \ i=1,2,\ldots, k.
    \]
    Consider indices $m_1,\ldots, m_d$ such that the arcs $T(\beta_{m_1}\otimes E_{m_1}),\ldots$, $T(\beta_{m_d}\otimes E_{m_d})$ form a cycle of length $d$ in $G(T(A))$.

    Let $T(E_i) = \alpha_i\otimes E_i'$ for all $i=1,2,\ldots, k$. Consider the matrix
    \[
    B = \left(\bigoplus\limits_{i\in \{m_1,\ldots, m_d\}} \alpha_{i}^{\otimes -1}\otimes  E_i\right)\oplus\left(\bigoplus\limits_{i\notin \{m_1,\ldots, m_d\}} \alpha_{i}^{\otimes -1}\otimes(-1)\otimes E_i\right).
    \]
    Note that the graph $G(T(B))$ contains a cycle of length $d$ with average weight 0, while the average weight of all other cycles in this graph is less than 0. Hence the critical subgraph $G_c(T(B))$ is a cycle of length $d$, and $\sigma(T(B)) = d<k=\sigma(B)$. This contradicts the fact that $T$ preserves the value $k$ of the cyclicity index.

    Thus, the graph $G(T(A))$ cannot contain cycles of length less than $k$. At the same time, $G(T(A))$ has exactly $k$ arcs in total, and this graph contains at least one cycle. Hence $G(T(A))$ is a cycle of length~$k$, and the mapping $\widehat T$ preserves the value $k$ of cycle lengths.
\end{proof}

\begin{corollary}\label{cor_adj_edges_k}
        Let $T: M_n\to M_n$ be a bijective linear mapping of tropical matrices preserving the value $k$ of the cyclicity index, $k\leq n$. Then 
        \begin{enumerate}
            \item for any matrix $A\in\mathcal{N}_n$ we have $T(A)\in\mathcal{N}_n$;
            \item the mappings $T|_{\mathcal{N}_n}$, $T|_{\mathcal{D}_n}$, $\widehat T|_{\mathcal{\widehat N}_n}$, $\widehat T|_{\mathcal{\widehat D}_n}$ are well defined and bijective.
        \end{enumerate}
\end{corollary}

\begin{proof}
    For $k=1$ the statement follows from Proposition~\ref{prop_diagonal_to_diagonal} and Corollary~\ref{bijective_on_diagonal_cells_cor}.

    Now let $k>1$. Let us prove the first part of the statement. Consider an arbitrary matrix $A\in\mathcal{N}_n$. By Lemma~\ref{surjective_maps}, $T(A) \in \mathcal{E}_n$. It is clear that there exist matrices $E_1, E_2,\ldots, E_{k-1}\in\mathcal{N}_n$ such that the graph $G(A\oplus E_1\oplus\ldots\oplus E_{k-1})$ is a cycle of length~$k$. By Lemma~\ref{lem_preserving_k_cycles}, the graph $G(T(A\oplus E_1\oplus\ldots\oplus E_{k-1}))$ is also a cycle of length~$k$, so it contains no loops. Thus $T(A)\in\mathcal{N}_n$, and statement~$(1)$ is proved.

    Statement $(2)$ follows from $(1)$ similarly to the proof of Corollary~\ref{bijective_on_diagonal_cells_cor}.
\end{proof}

\begin{corollary}\label{cor_permutation_indices_preserving_k}
    Let $T: M_n\to M_n$ be a bijective linear mapping of tropical matrices preserving the value $k\geq3$ of the cyclicity index, $n\geq\max(k, 4)$. Then there exists a permutation $\tau\in S_n$ such that exactly one of the following two conditions holds simultaneously for all $1\le i\ne j\le n$:
    \[
    T(E_{ij})=\alpha_{ij}\otimes E_{\tau(i)\tau(j)} \quad \text{or} \quad T(E_{ij})=\alpha_{ij}\otimes E_{\tau(j)\tau(i)}, \quad \alpha_{ij} \in \mathbb{R}.
    \]
\end{corollary}

\begin{proof}
    The desired statement follows immediately from Lemma~\ref{lem_preserving_k_cycles}, Corollary~\ref{cor_adj_edges_k}, and Theorem~\ref{th_vertices_permutation}.
\end{proof}

Recall that, for an arbitrary bijective linear mapping~$T: M_n\to M_n$, the numbers $\alpha_{ij}\in\mathbb{R}$ were defined in Eq.~\eqref{equation:basis-images}. Namely, for any $1\leq i,j\leq n$ there exist indices $k = k(i,j)$ and $m=m(i,j)$, $1\leq k,m\leq n$, such that $T(E_{ij}) = \alpha_{ij}\otimes E_{km}$.

\begin{notation}
    Let $C\in M_n$ be a matrix such that $G(C)$ and $G(T(C))$ are cycles of the same length $l\geq 2$. Then we denote 
    \[
    \Delta(C) = a(G(T(C))) - a(G(C)).
    \]
\end{notation}

\begin{lemma}\label{formula_delta_C}
    Let $T: M_n\to M_n$ be a bijective linear mapping of tropical matrices, and let $C\in M_n$ be a matrix such that the graphs $G(C)$ and $G(T(C))$ are cycles of length $l\geq 2$. Then
    \[
    \Delta(C) = \frac{1}{l}\sum_{\substack{(i,j)\\ e_{ij}\in G(C)}} \alpha_{ij}.
    \]

    In particular, if $C'\in M_n$ is a matrix such that $\widehat G(C') = \widehat G(C)$, then $\Delta(C')=\Delta(C)$.
\end{lemma}

\begin{proof}
    By Lemma~\ref{surjective_maps}, the mapping $T$ permutes the matrices from $\mathcal{E}_n$. Let $c_{ij}$ be the weight of the arc $e_{ij}$ in the graph $G(C)$. Then the image of this arc in $G(T(C))$ has weight $c_{ij}+\alpha_{ij}$. Therefore,
    \[
    a(G(T(C))) - a(G(C)) = \frac{1}{l}\sum_{e_{ij}\in G(C)} \big(c_{ij}+\alpha_{ij}-c_{ij}\big) = \frac{1}{l}\sum_{e_{ij}\in G(C)} \alpha_{ij}.
    \]

    To prove the second part of the lemma, it suffices to note that the right-hand side of the equality does not depend on the weights of the arcs in the graph $G(C)$.
\end{proof}

\begin{proposition}\label{lem_the_same_add_to_average_weight_k_prop7}
    Let $T: M_n\rightarrow M_n$ be a bijective linear transformation of tropical matrices and $\lambda\in\mathbb{R}$. Suppose that for all matrices $C\in M_n$ such that $G(C)$ is a cycle of length $2$ or $3$, the graph $G(T(C))$ is a cycle of the same length, and $\Delta(C) = \lambda$. Then there exist $u_1,\ldots, u_n \in \mathbb{R}$ such that $\alpha_{ij} = \lambda + u_j - u_i$ for all $1\leq i\neq j \leq n$.
\end{proposition}

\begin{proof} 
    For any pairwise distinct $1\leq i,j,l\leq n$ consider the cycles $i \to j \to i$ and $i \to j \to l \to i$. Then from the hypothesis and Lemma~\ref{formula_delta_C} we obtain
    \[
    \frac{\alpha_{ij}+\alpha_{ji}}{2} =
    \frac{\alpha_{ij}+\alpha_{jl}+\alpha_{li}}{3} = \lambda.
    \]

    Set $\alpha_{ij}' = \alpha_{ij} - \lambda$. Then 
    \[
    \alpha_{ij}'+\alpha_{ji}' = \alpha_{ij}'+\alpha_{jl}'+\alpha_{li}' = 0.
    \]
    
    Define $u_1 = 0$ and $u_i = \alpha_{1i}'$ for all $i=2,\ldots, n$. If $1, i, j$ are pairwise distinct, then $\alpha_{1j}' + \alpha_{j1}' = \alpha_{1i}' + \alpha_{ij}' + \alpha_{j1}' = 0$. Hence
    \[
    \alpha_{ij} = \lambda + \alpha_{ij}' = \lambda - \alpha_{1i}' - \alpha_{j1}' = \lambda + \alpha_{1j}' - \alpha_{1i}' = \lambda + u_j - u_i. 
    \]
    If $i=1$, then
    \[
    \alpha_{1j} = \lambda + \alpha_{1j}' = \lambda + u_j - u_1.
    \]
    If $j=1$, then
    \[
    \alpha_{i1} = \lambda + \alpha_{i1}' = \lambda - \alpha_{1i}' = \lambda + u_1 - u_i. \tag*{\qedhere}
    \]
\end{proof}

The following statements are used for the proof of Lemma~\ref{lem_the_same_add_to_average_weight_k} which is an analogue of Lemma~\ref{lem_the_same_add_to_average_weight}. In particular, we show that if $T$ preserves the value $k\geq 3$ of the cyclicity index and $n\geq\max(k, 4)$, then $T$ satisfies the hypotheses of Proposition~\ref{lem_the_same_add_to_average_weight_k_prop7}.

\begin{proposition}\label{lem_the_same_add_to_average_weight_k_prop3}
    Let $T: M_n\rightarrow M_n$ be a bijective linear transformation of tropical matrices preserving the value $k\geq 3$ of the cyclicity index, $n\geq\max(k, 4)$, and let matrices $C_1, C_2\in M_n$ be such that the graphs $\widehat G(C_1), \widehat G(C_2)$ are cycles of length $k$ and $2$, respectively, with a common arc. Then $\Delta(C_1) = \Delta(C_2)$.
\end{proposition}

\begin{proof}
    First we prove that $\Delta(C_1)\geq \Delta(C_2)$. Assume the contrary, and let $\mu = \Delta(C_2) - \Delta(C_1) > 0$. Consider matrices $C_1', C_2'\in M_n$ whose unweighted graphs $\widehat G(C_1'), \widehat G(C_2')$ coincide with $\widehat G(C_1)$ and $\widehat G(C_2)$, respectively, while the weights of all arcs of $G(C_1')$ and the weight of the arc of $G(C_2')$ which is common with $G(C_1')$ are equal to 0, and the arc of $G(C_2')$ which is not common with $G(C_1')$ has weight $-2\mu$. Then, by Lemma~\ref{formula_delta_C}, we have $\Delta(C_1') = \Delta(C_1)$ and $\Delta(C_2')=\Delta(C_2)$. Moreover,
    \[a(G(C_1')) = 0 > -\frac{2\mu}{2} = -\mu = a(G(C_2')).\]
    The graph $G(C_1'\oplus C_2')$ contains two cycles whose lengths are $k$ and $2$. Hence its critical subgraph is $G_c(C_1'\oplus C_2') = G(C_1')$, and $\sigma(C_1'\oplus C_2') = k$. At the same time, 
    \begin{align*}
        a(G(T(C_2'))) &= a(G(C_2')) + \Delta(C_2') = -\mu+\Delta(C_2) \\
        &= -\mu+\mu+\Delta(C_1) = \Delta(C_1') = a(G(T(C_1'))).
    \end{align*}

    By Corollary~\ref{cor_permutation_indices_preserving_k}, the graph $G(T(C_1'\oplus C_2'))$ contains exactly two cycles whose lengths are $k$ and $2$, with a common arc. Hence both cycles $G(T(C_1'))$ and $G(T(C_2'))$ are critical in this graph. Then $\sigma(T(C_1'\oplus C_2'))$ is a divisor of the length of the cycle $G(T(C_2'))$, that is, it divides 2, and, therefore, $\sigma(T(C_1'\oplus C_2'))\neq k$. This contradicts the fact that $T$ preserves the value $k$ of the cyclicity index. Hence $\Delta(C_1)\geq \Delta(C_2)$.

    Now we prove the desired equality. Consider the matrix $X=C_1\oplus C_1^t$. The following equality describes how the sum of the weights of all arcs in $G(X)$ is changed under the mapping $T$:
    \[
    w(G(T(X))) - w(G(X)) = k\cdot\Delta(C_1)+k\cdot\Delta(C_1^t).
    \]

    On the other hand, the graph $G(X)$ can be partitioned into $k$ cycles of length~$2$. Denote them by $C_{2,1}, C_{2,2}, C_{2,3}, \ldots, C_{2,k}$. By what has been proved, we obtain the inequality $\Delta(C_{2,i})\leq \min(\Delta(C_1),\Delta(C^t_1))$ for all $i=1,2,\ldots, k$. Then
    \begin{align*}
    w(G(T(X))) - w(G(X)) = \sum_{i=1}^{k} 2\Delta(C_{2,i}) &\leq 2k\cdot\min(\Delta(C_1),\Delta(C^t_1))\\
    &\leq w(G(T(X))) - w(G(X)).
    \end{align*}
    Hence $\Delta(C_{2,i})= \Delta(C_1) = \Delta(C^t_1)$ for all $i=1,2,\ldots, k$.
\end{proof}

\begin{proposition}\label{lem_the_same_add_to_average_weight_k_prop4}
    Let $T: M_n\rightarrow M_n$ be a bijective linear transformation of tropical matrices preserving the value $k\geq 3$ of the cyclicity index, $n\geq\max(k, 4)$. Then
    \begin{enumerate}
        \item for any matrices $C_1, C_2\in M_n$ such that $\widehat G(C_1), \widehat G(C_2)$ are cycles of length $2$, we have $\Delta(C_1) = \Delta(C_2)$.
        \item for any matrices $C_1, C_2\in M_n$ such that $\widehat G(C_1), \widehat G(C_2)$ are cycles of length $k$, we have $\Delta(C_1) = \Delta(C_2)$.
     \end{enumerate}
\end{proposition}

\begin{proof}
    (1) First, let the cycles $\widehat G(C_1), \widehat G(C_2)$ of length~$2$ have a common vertex. Consider a matrix $C\in M_n$ such that $\widehat G(C)$ is a cycle of length $k$ which has a common arc with each of these cycles. By Proposition~\ref{lem_the_same_add_to_average_weight_k_prop3}, $\Delta(C_1) = \Delta(C) = \Delta(C_2)$.

    Now let the cycles $\widehat G(C_1)$ and $\widehat G(C_2)$ have no common vertices. Then there exist pairwise distinct indices $1\leq i,j,m,l\leq n$ such that $\widehat G(C_1)$ is the cycle $i\to j\to i$, and $\widehat G(C_2)$ is the cycle $m\to l\to m$. Consider the matrix $C_3 = E_{im}\oplus E_{mi}$. By what has been proved in the previous paragraph, $\Delta(C_1) = \Delta(C_3) = \Delta(C_2)$.

    (2) Now let $\widehat G(C_1)$ and $\widehat G(C_2)$ be cycles of length $k$. Consider $C_1',C_2'\in M_n$ such that $\widehat G(C_i')$ is a cycle of length 2, and the graphs $\widehat G(C_i')$ and $\widehat G(C_i)$ have a common arc for each index $i=1,2$. By statement~$(1)$, we have $\Delta(C_1') = \Delta(C_2')$. Moreover, by Proposition~\ref{lem_the_same_add_to_average_weight_k_prop3}, we obtain $\Delta(C_1') = \Delta(C_1)$ and $\Delta(C_2') = \Delta(C_2)$. Hence the equality $\Delta(C_1) = \Delta(C_2)$ holds.
\end{proof}

\begin{proposition}\label{lem_the_same_add_to_average_weight_k_prop5}
    Let $T: M_n\rightarrow M_n$ be a bijective linear transformation of tropical matrices preserving the value $k\geq 3$ of the cyclicity index, $n\geq\max(k, 4)$. Let also matrices $C_1, C_2\in M_n$ be such that the graphs $\widehat G(C_1), \widehat G(C_2)$ are cycles of length $k$ and~$3$, respectively, at least one arc of the cycle $\widehat G(C_2)$ is common with $\widehat G(C_1)$, and the third vertex of $\widehat G(C_2)$ also belongs to $\widehat G(C_1)$. Then $\Delta(C_1) = \Delta(C_2)$.
\end{proposition}

\begin{proof}
    1. For $k=3$ the statement follows from Proposition~\ref{lem_the_same_add_to_average_weight_k_prop4}. In what follows we assume that $k\geq 4$. Since one arc of the cycle $\widehat G(C_2)$ is common with $\widehat G(C_1)$, and the third vertex also belongs to the cycle $\widehat G(C_1)$, without loss of generality we may assume that $\widehat G(C_1)$ is the cycle $v_1\to v_2\to\ldots\to v_k \to v_1$, and $\widehat G(C_2)$ is the cycle $v_1\to v_2\to v_i\to v_1$ for some $3\leq i\leq k$.

    2. First we prove that $\Delta(C_1)\geq\Delta(C_2)$. Assume the contrary, and denote $\mu = \Delta(C_2) - \Delta(C_1) > 0$. Let first $3<i<k$. Consider matrices $C_1', C_2'\in M_n$ such that $\widehat G(C_1') = \widehat G(C_1)$, $\widehat G(C_2') = \widehat G(C_2)$, the weights of all arcs of the cycle $G(C_1')$ and the weight of the arc $v_1 \to v_2$ of the cycle $G(C_2')$ are equal to $0$, and the arcs $v_2\to v_i$ and $v_i\to v_1$ of the cycle $G(C_2')$ have weight $-3\mu/2$. The graph $G(C_1'\oplus C_2')$ contains four cycles whose lengths are $k, 3, i, k-i+3$.  By Corollary~\ref{cor_permutation_indices_preserving_k}, the graph $G(T(C_1'\oplus C_2'))$ contains four cycles of the same lengths with the same common arcs. Moreover, $a(G(C_1')) = 0$ and
    \[
    a(G(C_2')) = \frac{0 -\frac{3\mu}{2}-\frac{3\mu}{2}}{3} = -\mu < 0.
    \]
    The average weights of the cycles of lengths $i$ and $k-i+3$ are also less than~$0$. Thus, $\sigma(C_1'\oplus C_2') = k$.

    By Lemma~\ref{formula_delta_C}, we obtain $\Delta(C_1') = \Delta(C_1)$ and $\Delta(C_2')=\Delta(C_2)$. In $G(T(C_1'\oplus C_2'))$ we have 
    \begin{align*}
    a(G(T(C_2'))) &= a(G(C_2')) + \Delta(C_2') = -\mu + \Delta(C_2) \\
    &= -\mu + \mu + \Delta(C_1) = \Delta(C_1') = a(G(T(C_1'))).
    \end{align*}

    Therefore, the cycles $G(T(C_2'))$ and $G(T(C_1'))$ are either both critical or both non-critical in the graph $G(T(C_1'\oplus C_2'))$. In the first case $\sigma(T(C_1'\oplus C_2'))$ divides $3$. Hence we have $\sigma(T(C_1'\oplus C_2'))\neq k$.

    In the second case the critical subgraph $G_c(T(C_1'\oplus C_2'))$ can consist either only of the cycle of length $i$, or only of the cycle of length $k-i+3$, or of both these cycles at once. The first two possibilities imply that the cyclicity index of the matrix $T(C_1'\oplus C_2')$ is equal to $i<k$ or $k-i+3<k$, respectively. The third one means that the graph $G_c(T(C_1'\oplus C_2'))$ consists of the cycles $G(T(C_2'))$ and $G(T(C_1'))$, and, similarly to the first case, we obtain $\sigma(T(C_1'\oplus C_2'))\neq k$. We arrive at a contradiction with the fact that the mapping $T$ preserves the value $k$ of the cyclicity index.

    3. It remains to consider the cases $i=3$ and $i=k$. We assume that $i=3$; the second case is treated similarly. Consider matrices $C_1'', C_2''\in M_n$ such that $\widehat G(C_1'')~=~\widehat G(C_1)$, $\widehat G(C_2'') = \widehat G(C_2)$, the weights of all arcs of the cycle $G(C_1'')$ and the weights of the arcs $v_1 \to v_2$ and $v_2 \to v_3$ of the cycle $G(C_2'')$ are equal to $0$, and the arc $v_3\to v_1$ of the cycle $G(C_2'')$ has weight $-3\mu$. The graph $G(C_1''\oplus C_2'')$ consists of two cycles whose lengths are $k$ and~$3$.  By Corollary~\ref{cor_permutation_indices_preserving_k}, the graph $G(T(C_1''\oplus C_2''))$ also contains exactly two cycles whose lengths are $k$ and $3$, with a common arc. We have $a(G(C_1'')) = 0$ and 
    \[
    a(G(C_2'')) = \frac{0+0-3\mu}{3} = -\mu <0.
    \]
    Hence $\sigma(C_1''\oplus C_2'') = k$. 

    From Lemma~\ref{formula_delta_C} we obtain $\Delta(C_1'') = \Delta(C_1)$ and $\Delta(C_2'')=\Delta(C_2)$. Similarly to the case $3<i<k$, in $G(T(C_1''\oplus C_2''))$ we have 
    \begin{align*}
    a(G(T(C_2''))) &= a(G(C_2'')) + \Delta(C_2'') = -\mu + \Delta(C_2) \\
    &= -\mu + \mu + \Delta(C_1) = \Delta(C_1'') = a(G(T(C_1''))),
    \end{align*}
    and both cycles $G(T(C_2''))$ and $G(T(C_1''))$ are critical in the graph $G(T(C_1''\oplus C_2''))$. Then $\sigma(T(C_1''\oplus C_2''))$ divides $3$, so $\sigma(T(C_1''\oplus C_2''))\neq k$. Thus, we again arrive at a contradiction with the fact that the mapping $T$ preserves the value~$k$ of the cyclicity index. The inequality $\Delta(C_1)\geq\Delta(C_2)$ is proved.

    4. We partition the cycle of length $k$ into $k-2$ cycles of length 3:
    \[
    v_1\to v_2\to v_3\to v_1,\quad v_1\to v_3\to v_4\to v_1,\quad
    \ldots,\quad v_1\to v_{k-1}\to v_k\to v_1.
    \]

    Denote the tropical adjacency matrices of these cycles by $C_{3,3},\ldots, C_{3,k}$. Besides, denote by $C_{2,3}, \ldots, C_{2, k-1}$ the tropical adjacency matrices of the cycles $v_1\to v_3\to v_1,$ $\ldots, v_1\to v_{k-1}\to v_1$. By Propositions~\ref{lem_the_same_add_to_average_weight_k_prop3} and~\ref{lem_the_same_add_to_average_weight_k_prop4}, we obtain $\Delta(C_{2,i}) = \Delta(C_1)$ for any $i=3,\ldots, k-1$. Then 
    \begin{align*}
    \sum_{i=3}^{k}3\Delta(C_{3,i}) &= k\cdot\Delta(C_1) + \sum_{i=3}^{k-1}2\Delta(C_{2,i}) \\
    &= k\cdot\Delta(C_1) + 2(k-3)\cdot\Delta(C_1) = 3(k-2)\Delta(C_1).
    \end{align*}

    At the same time, by Steps 2--3, we have $\Delta(C_1)\geq\Delta(C_{3,i})$ for all $i=3,\ldots, k$. Thus the equality above is possible only if $\Delta(C_1)=\Delta(C_{3,i})$ for all $i=3,\ldots, k$. 
    
    5. The argument of Step 4 remains valid for an arbitrary choice of the initial vertex of $C_1$. In order to obtain the desired equality for the original triangle $v_1\to v_2\to v_i\to v_1$, we cyclically relabel the vertices:
    \[
    (w_1,\ldots,w_k) =(v_i,\ldots,v_k,v_1,\ldots,v_{i-1}).
    \]
    Then this triangle has the form
    \[
    w_1\to w_j\to w_{j+1}\to w_1,
    \qquad j=k-i+2,\quad 2\le j\le k-1.
    \]
    Applying the result of Step 4 to this numbering of the vertices, we get $\Delta(C_2)=\Delta(C_1)$.
\end{proof}

\begin{proposition}\label{lem_the_same_add_to_average_weight_k_prop6}
    Let $T: M_n\rightarrow M_n$ be a bijective linear transformation of tropical matrices preserving the value $k\geq 3$ of the cyclicity index, $n\geq\max(k, 4)$, and let matrices $C_1, C_2\in M_n$ be such that the graphs $\widehat G(C_1), \widehat G(C_2)$ are cycles of length $3$. Then $\Delta(C_1) = \Delta(C_2)$.
\end{proposition}

\begin{proof}
    The statement follows immediately from Propositions~\ref{lem_the_same_add_to_average_weight_k_prop4}(2) and~\ref{lem_the_same_add_to_average_weight_k_prop5}.
\end{proof}

Finally, we are ready to prove the following lemma.

\begin{lemma}\label{lem_the_same_add_to_average_weight_k}
    Let $T: M_n\rightarrow M_n$ be a bijective linear transformation of tropical matrices preserving the value $k\geq 3$ of the cyclicity index, $n\geq\max(k, 4)$. Then there exists a number $\lambda\in\mathbb{R}$ such that $\Delta(C) = \lambda$ for all matrices $C\in M_n$ whose graph $\widehat G(C)$ is a cycle of length at least $2$. Equivalently, $a(G(T(C))) = \lambda\otimes a(G(C))$.
\end{lemma}

\begin{proof}
    From Proposition~\ref{lem_the_same_add_to_average_weight_k_prop4}(2) it follows that one can define a number $\lambda\in\mathbb{R}$ such that $\Delta(C_k) = \lambda$ for any matrix $C_k\in M_n$ whose graph $\widehat G(C_k)$ is a cycle of length $k$. Propositions~\ref{lem_the_same_add_to_average_weight_k_prop3}--\ref{lem_the_same_add_to_average_weight_k_prop6} imply that $\Delta(C_{2,3}) = \lambda$ for all matrices $C_{2,3}\in M_n$ whose graph $\widehat G(C_{2,3})$ is a cycle of length $2$ or $3$.

    Then, by Proposition~\ref{lem_the_same_add_to_average_weight_k_prop7}, there exist $u_1, \ldots, u_n \in \mathbb{R}$ such that $\alpha_{ij} = \lambda + u_j - u_i$ for all $1 \leq i \neq j \leq n$.

    Hence, for any matrix $C\in M_n$ whose graph $\widehat G(C)$ is a cycle of length at least 2, by Lemma~\ref{formula_delta_C} we have
    \begin{align*}
    \Delta(C) &= a(G(T(C))) - a(G(C)) = \frac{1}{|G(C)|}\sum_{e_{ij}\in G(C)} \alpha_{ij}\\
    &= \frac{1}{|G(C)|}\sum_{e_{ij}\in G(C)} (\lambda + u_j - u_i) \\
    &= \frac{1}{|G(C)|}\big(\lambda\cdot |G(C)| + \sum_{e_{ij}\in G(C)} (u_j - u_i)\big) = \lambda. \tag*{\qedhere}
    \end{align*}
\end{proof}

\begin{theorem}\label{th_non_diagonal_elements_k}
Let $T : M_n \to M_n$ be a bijective linear transformation
of tropical matrices preserving the value $k \ge 3$ of the cyclicity
index, $n \ge \max(k,4)$. Then there exist
\begin{itemize}
    \item a matrix $B = (b_{ij})\in M_n$ such that $b_{ij}\neq-\infty$ for all $1\leq i,j\leq n$, and all cycles in $G(B)$ have weight $0$,
    \item a number $\lambda\in\mathbb{R}$,
    \item a permutation $\tau \in S_n$,
\end{itemize}
such that exactly one of the following two conditions holds simultaneously for all $1\le i\ne j\le n$:
\[
T(E_{ij})=\lambda \otimes b_{ij}\otimes E_{\tau(i)\tau(j)} \quad \text{or} \quad T(E_{ij})=\lambda \otimes b_{ij}\otimes E_{\tau(j)\tau(i)}.
\]
\end{theorem}

\begin{proof}
By Corollary~\ref{cor_permutation_indices_preserving_k}, there exists a permutation $\tau\in S_n$ such that one of the equalities
\[
    T(E_{ij})=\alpha_{ij}\otimes E_{\tau(i)\tau(j)}
    \quad \text{or} \quad T(E_{ij})=\alpha_{ij}\otimes E_{\tau(j)\tau(i)}
\]
holds simultaneously for all $1\le i\ne j\le n$.

By Lemma~\ref{lem_the_same_add_to_average_weight_k}, there exists a number $\lambda\in\mathbb{R}$ such that for
any matrix $C\in M_n$ whose graph $\widehat G(C)$ is a cycle of length at least $2$, we have
$a(G(T(C))) = \lambda \otimes a(G(C))$.
Set
\[
b_{ij} = 
\begin{cases}
    \alpha_{ij}\otimes\lambda^{\otimes-1}, & i \neq j,\\
    0, & i = j.
\end{cases}
\]
Then $b_{ij}\ne -\infty$ for all $1\le i,j\le n$ and $\alpha_{ij} = \lambda\otimes b_{ij}$ for $i\neq j$. Similarly to the proof of Theorem~\ref{main_theorem}, we obtain that all cycles in $G(B)$ have weight $0$. 
\end{proof}

\subsection{Action on the diagonal entries of matrices}

\begin{notation}
    Let 
    \begin{equation} \label{equation:primes}
    k = \prod\limits_{\nu=1}^{r} p_\nu^{\alpha_\nu},
    \end{equation}
    where $p_1, p_2,\ldots, p_r$ are distinct primes and $\alpha_{\nu} \in \mathbb{N}$ for all $1 \leq \nu \leq r$. We denote:
    \[
    m(k) = \sum_{\nu=1}^{r} p_\nu^{\alpha_\nu} \quad \text{and} \quad
    q(k) = \min\{d\in\mathbb{N}: d\nmid k\}.
    \]
\end{notation}

\begin{lemma}\label{lem_m(k)}
    Let $d_1, d_2, \ldots, d_n \in \mathbb{N}$ be such that their least common multiple satisfies $\lcm(d_1,\ldots, d_n) = k > 1$. Then
    \[
    d_1 + d_2 + \ldots + d_n \geq m(k),
    \]
    and equality holds if and only if $n = r$ and there exists a permutation $\tau \in S_r$ such that $d_{\tau(\nu)} = p_\nu^{\alpha_\nu}$ for all $1 \leq \nu \leq r$.
\end{lemma}

\begin{proof}
    It follows from the equality $\lcm(d_1,\ldots, d_n) = k$ that for each $\nu$ in factorization~\eqref{equation:primes} there exists an index $j_\nu$ such that $p_\nu^{\alpha_\nu}\mid d_{j_\nu}$. For each $1 \leq j \leq n$ we set
    \[
    I_j = \{\nu: j_\nu = j\}.
    \]
    Let us prove that for all $1\leq j\leq n$ we have
    \[
    d_j\geq\sum_{\nu\in I_j} p_\nu^{\alpha_\nu}.
    \]
    If $I_j=\varnothing$, then the inequality is obvious, and it is strict. Now let $I_j\neq\varnothing$. The numbers $p_\nu^{\alpha_\nu}$, $\nu\in I_j$, are pairwise coprime, and each of them divides $d_j$. Hence
    \[
    \prod_{\nu\in I_j} p_\nu^{\alpha_\nu}\mid d_j, \quad \text{and thus} \quad
    d_j\geq\prod_{\nu\in I_j} p_\nu^{\alpha_\nu}.
    \]

    Note that for any integers $x_1, x_2\geq 2$ we have $x_1x_2\geq x_1 + x_2$, and equality is attained if and only if $x_1 = x_2 = 2$. By induction, for $t \in \mathbb{N}$ and $x_1,\ldots, x_t\geq 2$ we have
    \[
    \prod_{i=1}^t x_i \geq\sum_{i=1}^t x_i,
    \]
    and equality holds if and only if either $t = 1$, or $t = 2$ and $x_1 = x_2 = 2$. Since $p_\nu^{\alpha_\nu}\geq 2$, we obtain
    \begin{equation} \label{equation:primes-inequality}
        d_j\geq\prod_{\nu\in I_j} p_\nu^{\alpha_\nu}\geq\sum_{\nu\in I_j} p_\nu^{\alpha_\nu}.
    \end{equation}
    Since the primes $p_\nu$ are pairwise distinct, inequality~\eqref{equation:primes-inequality} becomes an equality only if $I_j$ consists of a single element~$\nu$ and $d_j = p_\nu^{\alpha_\nu}$. Summing these inequalities over all $j$, we obtain
    \[
    \sum_{j=1}^n d_j\geq \sum_{j=1}^n\sum_{\nu\in I_j} p_\nu^{\alpha_\nu} = \sum_{\nu=1}^r p_\nu^{\alpha_\nu} = m(k).
    \]
    If equality is attained, then there are no $d_j$ with $I_j = \varnothing$, and for each of the remaining~$d_j$ Eq.~\eqref{equation:primes-inequality} becomes an equality. Hence in this case $n = r$, and the permutation $\tau \in S_r$ can be defined by the formula $\tau(\nu) = j_\nu$.
\end{proof}

\begin{corollary}\label{cor_m(k)}
    Let $G$ be a directed graph with the cyclicity index $\sigma_G=k > 1$. Then~$G$ contains at least $m(k)$ vertices.
\end{corollary}

\begin{proof}
    Consider the strongly connected components of the graph~$G$ with cyclicity indices $d_1, \ldots, d_n > 1$. Their lcm equals $k$ by the definition of the cyclicity index. In the component whose index equals $d_j$, every cycle has length divisible by $d_j$. In particular, the number of vertices in this strongly connected component is at least $d_j$. Then the desired statement follows from Lemma~\ref{lem_m(k)}.
\end{proof}

Let $A\in M_n$ be an arbitrary tropical matrix with $\sigma(A) = k > 1$. Then the graph $G(A)$ contains a critical cycle of length at least $2$. We denote the largest average weight of a cycle in $G(A)$ by~$\mu$. It is clear that $A_{ii}\leq\mu$ for all $1\leq i\leq n$. 

Let  $H = H(A)$ be the union of all critical cycles of length at least $2$, and let
\[
L = L(A) = \{i: A_{ii}=\mu\}
\] 
be the set of vertices incident to critical loops. Then $G_c(A)$ is obtained from $H$ by adding loops at the vertices from $L$. Consider the strongly connected components $H_1, H_2,\ldots, H_s$ of the graph $H$ and set
\begin{equation}\label{equation:H_j-cyclicity-vertices}
d_j = \sigma_{H_j},\quad V_j = V(H_j).
\end{equation}

Note that the cyclicity index of the strongly connected component of $G_c(A)$ which corresponds to $H_j$ equals
\[
\begin{cases}
    1, & V_j\cap L\neq\varnothing,\\
    d_j, & V_j\cap L=\varnothing.
\end{cases}
\]
Hence
\begin{equation} \label{equation:cyclicity-loops}
    k = \sigma(A) = \lcm\{d_j: V_j\cap L=\varnothing\}.
\end{equation}

\begin{lemma}\label{lem_m(k)_2}
    Let $A\in M_n$, $\sigma(A) = k > 1$, and let $H, L, H_j, V_j, d_j$ be defined as above. Then the following statements hold:
    \begin{enumerate}
        \item The union of the components $H_j$ containing no vertices from $L$ contains at least $m(k)$ vertices.

        \item If $n = m(k)$, then $L = \varnothing$.

        \item If $n < m(k) + q(k)$, then $d_j \mid k$ for all $1\leq j\leq s$.
    \end{enumerate}
\end{lemma}

\begin{proof}
\begin{enumerate}
    \item Follows from Eq.~\eqref{equation:cyclicity-loops} and Corollary~\ref{cor_m(k)} applied to the union of the components $H_j$.

    \item If $n = m(k)$, then by item~(1) we obtain that every vertex  belongs to some component $H_j$ containing no vertices from $L$. Hence $L = \varnothing$.

    \item Now suppose that $n < m(k) + q(k)$ and there exists an index $t$ such that $d_t\nmid k$. By Eq.~\eqref{equation:cyclicity-loops}, 
    we obtain that $H_t$ contains a vertex from $L$. Moreover,
    \[
    |V_t|\geq d_t\geq q(k).
    \]

    Thus,
    \[
    |V(G(A))|\geq |V_t| + \sum_{j:V_j\cap L=\varnothing} |V_j|\geq q(k) + m(k),
    \]
    a contradiction. \qedhere
    \end{enumerate}
\end{proof}

\begin{theorem}\label{th_diag_k}
    Let $n \ge k\geq 3$, and let $R : M_n \to M_n$ be a bijective linear transformation of tropical matrices defined as follows:
    \[
R(E_{ij}) = \begin{cases}
    E_{ij},& 1\le i\ne j\le n,\\
    \lambda_{i}\otimes E_{\pi(i)\pi(i)},& 1\leq i=j\leq n,
\end{cases}
\]
where $\lambda_{i}\in\mathbb{R}$ and $\pi\in S_n$. Then $R$ preserves the value $k$ of the cyclicity index if and only if the condition for the corresponding case holds:
\begin{enumerate}[font=\upshape,label=(\roman*)]
    \item if $n = m(k)$, then $\lambda_i\leq 0$ for all $1\leq i\leq n$, and the permutation $\pi$ is arbitrary;
    \item if $m(k)<n<m(k)+q(k)$, then $\lambda_i\leq 0$ for all $1\leq i\leq n$ and, if $\lambda_i = 0$ for some index $i$, then $\pi(i) = i$;
    \item if $n\geq m(k) + q(k)$, then $\lambda_i= 0$ and $\pi(i) = i$ for all $1\leq i\leq n$. In this case $R$ is the identity mapping on the whole set $M_n$.
\end{enumerate}
\end{theorem}

\begin{proof}
    1. Suppose that the mapping $R$ preserves the value $k$ of the cyclicity index, and let us show that it satisfies the required conditions. First we prove that $\lambda_i\leq 0$ for all $1\leq i\leq n$.

    Assume the contrary: let there exist an index $i$ such that $\lambda_i > 0$. Consider a matrix $C_k\in M_n$ whose graph is a cycle of length $k$ with zero arc weights, and set
    \[
    A = C_k\oplus \left(-\frac{\lambda_i}{2} \right)\otimes E_{ii}.
    \]
    We have $\sigma(A)=k$. At the same time,
    \begin{align*}
    R(A) &= R(C_k)\oplus \left(-\frac{\lambda_i}{2}\right)\otimes R(E_{ii}) \\
    &= C_k\oplus \left(-\frac{\lambda_i}{2}\right)\otimes\lambda_i\otimes E_{\pi(i)\pi(i)} = C_k\oplus \frac{\lambda_i}{2}\otimes E_{\pi(i)\pi(i)}.
    \end{align*}
    Then $\sigma(R(A)) = 1$, hence the mapping $R$ does not preserve the value $k$ of the cyclicity index, a contradiction.

    2. Now let $n > m(k)$. Suppose that for some $i$ we have $\lambda_i = 0$ but $\pi(i)\neq i$. Let $k = \prod\limits_{\nu=1}^{r} p_\nu^{\alpha_\nu}$, where $p_1, \ldots, p_r$ are distinct primes. We choose $m(k)$ vertices different from the vertex $i$, among which there is the vertex $\pi(i)$. This is possible since $n > m(k)$ and $i\neq\pi(i)$. Consider matrices $C_1, \ldots, C_r$ such that their graphs are pairwise disjoint cycles of lengths $p_1^{\alpha_1}, \ldots, p_r^{\alpha_r}$ built on the chosen vertices, with zero arc weights. Without loss of generality we may assume that the vertex $\pi(i)$ belongs to the cycle $G(C_1)$. Set
    \[
    A = E_{ii}\oplus\bigoplus_{j=1}^r C_j.
    \]

    Then $G(A)$ consists of disjoint cycles of lengths $1,$ $ p_1^{\alpha_1},\ldots,p_r^{\alpha_r}$, each of weight~$0$. Hence
    \[
    \sigma(A) = \lcm(1, p_1^{\alpha_1}, \ldots,p_r^{\alpha_r}) = p_1^{\alpha_1} \ldots p_r^{\alpha_r} = k.
    \]

    Moreover,
    \[
    R(A) = R(E_{ii})\oplus\bigoplus_{j=1}^r R(C_j) = E_{\pi(i)\pi(i)}\oplus\bigoplus_{j=1}^r C_j.
    \]

    Since the vertex $\pi(i)$ belongs to the cycle $G(C_1)$, we have
    \[
     \sigma(R(A)) = \lcm(1, p_2^{\alpha_2},\ldots,p_r^{\alpha_r}) = p_2^{\alpha_2} \ldots p_r^{\alpha_r} = \frac{k}{p_1^{\alpha_1}} \neq k.
    \]

    Thus the mapping $R$ does not preserve the value $k$ of the cyclicity index, a contradiction. Hence, if $\lambda_i=0$, then $\pi(i) = i$.

    3. Next, suppose that $n\geq m(k)+q(k)$. Let us show that $\lambda_i = 0$ for all $1\leq i\leq n$. In Step 1 of this proof it was proved that $\lambda_i \leq 0$. Assume there exists an index $i$ such that $\lambda_i < 0$. 

    Consider matrices $C_1, C_2, \ldots, C_r, Q$ such that their graphs are pairwise disjoint cycles of lengths $p_1^{\alpha_1}, \ldots, p_r^{\alpha_r},$ $q(k)$, respectively, all of whose arcs have weight 0. We choose these matrices so that the vertex $i$ belongs to the cycle $G(Q)$. Set
    \[
    A = E_{ii}\oplus Q\oplus\bigoplus_{j=1}^r C_j.
    \]
    Note that the number of strongly connected components of the graph $G_c(A) = G(A)$ containing nontrivial cycles equals $r+1$, and their cyclicity indices are $1, p_1^{\alpha_1}, \ldots, p_r^{\alpha_r}$. All other strongly connected components, if any, correspond to isolated vertices and have cyclicity index $1$. Therefore
    \[
    \sigma(A) = \lcm(1, p_1^{\alpha_1}, \ldots, p_r^{\alpha_r}) = p_1^{\alpha_1}\ldots p_r^{\alpha_r} = k.
    \]

    Since the image of the loop $E_{ii}$ in $G(R(A))$ has weight $\lambda_i < 0$, it is no longer critical. Thus we obtain
    \begin{gather*}
    G_c(R(A)) = G\left(Q\oplus\bigoplus_{j=1}^r C_j\right),\\
    \sigma(R(A)) = \lcm(q(k), p_1^{\alpha_1}, \ldots, p_r^{\alpha_r}) = \lcm(q(k), k) > k.
    \end{gather*}

    The last inequality holds because $q(k)\nmid k$ by the definition of $q(k)$. This contradicts the fact that $R$ preserves the value $k$ of the cyclicity index. Hence, for such $n$, all $\lambda_i=0$. Then, by Step 2, $\pi = \mathrm{id}$. The necessity of conditions (i)--(iii) is proved.

    4. Now we prove the converse statement: if the mapping $R$ satisfies one of the conditions (i), (ii), or (iii), then it preserves the value $k$ of the cyclicity index.

    Consider an arbitrary matrix $A\in M_n$ with $\sigma(A) = k$. We introduce the following notation. Let $\mu=\mu(A)$ be the largest average weight in the graph $G(A)$, $H=H(A)$ be the union of all critical cycles of $G(A)$ of length at least $2$, and $L = L(A) = \{i: A_{ii}=\mu\}$ be the set of vertices incident to critical loops.

    In all three cases $\lambda_i\leq0$, so
    \[
    (R(A))_{\pi(i)\pi(i)} = A_{ii}+\lambda_i\leq\mu\quad\text{for all}\ 1\leq i\leq n.
    \] 
    Non-loop arcs do not change under the mapping $R$, therefore $\mu(R(A)) = \mu(A) = \mu$ and $H(R(A)) = H(A) = H$. Since  $A_{ii}\leq\mu$ and $\lambda_{i}\leq 0$, the condition $\pi(i)\in L' = L(R(A))$ implies $A_{ii}=\mu$ and $\lambda_{i}= 0$, that is,
    \[
    L' = \{\pi(i): i\in L, \; \lambda_i = 0\}.
    \]

    Next, we consider each case separately.

    \begin{enumerate}[font=\upshape,label=(\roman*)]
        \item If $n=m(k)$, then, by Lemma~\ref{lem_m(k)_2}(2), we have $L=\varnothing$. Hence $L'=\varnothing$, that is, the graphs $G(A)$ and $G(R(A))$ contain no critical loops. Then, by the definition of the mapping~$R$, the graphs $G_c(A)$ and $G_c(R(A))$ coincide for any permutation $\pi$, and thus $\sigma(A) = \sigma(R(A))$.

        \item By the assumption, $\lambda_i = 0$ implies $\pi(i) = i$. Therefore
        \[
        L' = \{\pi(i): i\in L, \; \lambda_i = 0\}= \{i\in L: \lambda_i=0\}\subseteq L.
        \]

        Thus, under the mapping $R$, some loops can disappear from the critical subgraph, but the rest of it remains unchanged. Recall that, by Eq.~\eqref{equation:cyclicity-loops},
        \[
        k = \sigma(A) = \lcm\{d_j: V_j\cap L=\varnothing\}
        \]
        where $d_j$ and $V_j$ are, respectively, the cyclicity index and the vertex set of the strongly connected component $H_j$ of the graph $H$, defined in~\eqref{equation:H_j-cyclicity-vertices}. Hence we obtain
        \[
        k\mid \sigma(R(A)) = \lcm\{d_j: V_j\cap L'=\varnothing\},
        \]
        since all strongly connected components $H_j$ whose lcm gave $k$ still participate in this lcm. On the other hand, by Lemma~\ref{lem_m(k)_2}(3), all indices $d_j$ divide $k$. Thus $\sigma(R(A)) \mid k$, whence $\sigma(R(A)) = k$.

        \item In this case $R$ is the identity mapping, and the statement is obvious. \qedhere
    \end{enumerate}
\end{proof}

\begin{remark} \label{remark:case-(i)}
    Applying Lemma~\ref{lem_m(k)} to $n = 1$ and $d_1 = k$, we obtain that $k \geq m(k)$, and equality is attained if and only if $k$ is a prime power. Since $n \geq k$, this condition is necessary for case~\textup{(i)} to be possible in Theorem~\ref{th_diag_k}. If $k$ is not a prime power, then only cases \textup{(ii)} and~\textup{(iii)} are possible.
\end{remark}

\begin{proposition} \label{proposition:case-(iii)}
    Suppose that $k > 1$ is not a prime power and $k \neq 6$. Then $k \geq m(k) + q(k)$. Therefore, for such~$k$, only case~\textup{(iii)} is possible in Theorem~\ref{th_diag_k}.
\end{proposition}

\begin{proof}
    By the assumption, $r \geq 2$ in factorization~\eqref{equation:primes}. Without loss of generality we may assume that $p_1 < p_2 < \ldots < p_r$.
    
    First suppose that at most one of the numbers $2$ and $3$ belongs to $p_1, \ldots, p_r$. Then $p_1 \geq 2$, $p_2 \geq 5$ and $q(k) \in \{ 2, 3\}$. Note that for $x_1 \geq 2$ and $x_2 \geq 5$ we have $x_1x_2 \geq x_1 + x_2 + 3$. Similarly to the proof of Lemma~\ref{lem_m(k)}, we obtain
    \[
    m(k) + q(k) \leq p_1^{\alpha_1} + p_2^{\alpha_2} + \ldots + p_r^{\alpha_r} + 3 \leq p_1^{\alpha_1} + p_2^{\alpha_2} \ldots p_r^{\alpha_r} + 3 \leq p_1^{\alpha_1} p_2^{\alpha_2}\ldots p_r^{\alpha_r} = k.
    \]

    Now let $p_1 = 2$ and $p_2 = 3$. If $k = 6$, then $m(k) = 2 + 3$ and $q(k) = 4$, therefore $m(k) + q(k) = 9 > k$. From now on we assume that $k \neq 6$. If $\alpha_1 = 1$, that is, $k$ is divisible by $2$ but not by $4$, then $q(k) = 4$ and $p_2^{\alpha_2} \ldots p_r^{\alpha_r} \geq 9$. Hence
    \[
    m(k) + q(k) = (2 + p_2^{\alpha_2} + \ldots + p_r^{\alpha_r}) + 4 \leq 6 + p_2^{\alpha_2} \ldots p_r^{\alpha_r} < 2 p_2^{\alpha_2} \ldots p_r^{\alpha_r} = k.
    \]
    
    Now consider the case $\alpha_1 \geq 2$. Let $p$ be the smallest prime greater than $p_r$. Then, as is well known, $p < 2p_r$, therefore $q(k) \leq p \leq 2p_r - 1$. Since $2^{\alpha_1} \geq 4$ and $p_r \geq 3$, we obtain $2p_r \leq 3(p_r - 1) \leq (2^{\alpha_1} - 1)(p_r^{\alpha_r} -1)$. Thus,
    \begin{align*}
        m(k) + q(k) &\leq (2^{\alpha_1} + p_2^{\alpha_2} + \ldots + p_r^{\alpha_r}) + (2p_r - 1) \\
    &\leq p_2^{\alpha_2} + \ldots + p_{r-1}^{\alpha_{r-1}} + 2^{\alpha_1}p_r^{\alpha_r} \leq 2^{\alpha_1} p_2^{\alpha_2} \ldots p_r^{\alpha_r} = k. \tag*{\qedhere}
    \end{align*}
\end{proof}

\begin{theorem}\label{th_final_k}
    Let $k\geq 3$, $n\geq\max(k, 4)$, and let $T:M_n\to M_n$ be a bijective linear
    transformation of tropical matrices. Then the following conditions are equivalent:
    \begin{enumerate}
        \item $T$ preserves the value $k$ of the cyclicity index;
        
        \item there exist
        \begin{itemize}
            \item a matrix $B = (b_{ij})\in M_n$ such that $b_{ij}\neq-\infty$ for all $1\leq i,j\leq n$, and all cycles in $G(B)$ have weight 0,
            \item a number $\lambda\in\mathbb{R}$,
            \item a diagonal matrix $\Lambda = \mathrm{diag}(\lambda_1,\ldots,\lambda_n)$, where $\lambda_1,\ldots,\lambda_n\in\mathbb{R}$,
            \item permutation matrices $P_1, P_2\in M_n$,
        \end{itemize}
    such that
        \[
        T = T_1\circ T_2\circ T_3\circ T_4\circ T_5\quad\text{or}\quad T = T_2\circ T_3\circ T_4\circ T_5,
        \]
        for the mappings $T_j$ defined by the formulas
        \begin{gather*}
            \begin{aligned}
            &T_1(A) = A^t, & T_2(A) = P_1\otimes A\otimes P_1^t,&\\
            &T_3(A) = \lambda\otimes A, & T_4(A) = A\circ B,&
            \end{aligned}\\
            T_5(A) = N(A)\oplus (P_2\otimes \Lambda\otimes  D(A)\otimes P_2^t)
        \end{gather*}
        for all $A\in M_n$, where $D(A)\oplus N(A) = A$, the matrix $D(A)$ is diagonal, and all diagonal entries of the matrix $N(A)$ are equal to $-\infty$.
        
        Moreover, the following restrictions hold:
        \begin{enumerate}[font=\upshape,label=(\roman*)]
            \item $\lambda_i\leq 0$ for all $1\leq i\leq n$,
            \item if $n > m(k)$ and $\lambda_i = 0$, then $(P_2)_{ii} = 0$,
            \item if $n\geq m(k)+q(k)$, then $\lambda_i = 0$ for all $1\leq i\leq n$; hence $P_2$ is the tropical identity matrix, and $T_5$ is the identity mapping.
        \end{enumerate}
    \end{enumerate}
\end{theorem}

\begin{proof}
    1. First, let the mapping $T$ preserve the value $k$ of the cyclicity index. By Theorem~\ref{th_non_diagonal_elements_k}, there exist a matrix~$B$, a number $\lambda$ and a permutation $\tau$ such that one of the equalities
    \[
    T(E_{ij})=\lambda \otimes b_{ij}\otimes E_{\tau(i)\tau(j)} \quad \text{or} \quad T(E_{ij})=\lambda \otimes b_{ij}\otimes E_{\tau(j)\tau(i)}
    \]
    holds simultaneously for all off-diagonal matrix units.
    
    Let $P_1$ be the matrix corresponding to the permutation $\tau^{-1}$. Then the values of $T$ on the off-diagonal matrix units coincide with the values of the mapping $T'$ which equals $T_2 \circ T_3 \circ T_4$ or $T_1 \circ T_2 \circ T_3 \circ T_4$, respectively. Note that all diagonal entries of the matrix~$B$ are equal to~$0$, therefore for all $1 \leq i \leq n$ we have
    \[
    T'(E_{ii}) = \lambda\otimes E_{\tau(i)\tau(i)}.
    \]
    In Theorem~\ref{main_theorem} and Corollary~\ref{cor_preserving_cyclicity} it was shown that the mappings $T_1, T_2, T_3, T_4$ are bijective and preserve the cyclicity index, and their inverse mappings also preserve the cyclicity index. Hence $T'$ and $(T')^{-1}$ preserve the cyclicity index.

    2. Set $R = (T')^{-1}\circ T$. Then $R$ preserves the value $k$ of the cyclicity index as a composition of two mappings preserving the value $k$ of the cyclicity index, and
    \[
        R(E_{ij}) = E_{ij}\quad\text{for all }i\neq j.
    \]
    By Corollary~\ref{cor_adj_edges_k}, the mappings $T$ and $T'$ act bijectively on the set $\mathcal{D}_n$. Then this is also true for $R$. Therefore there exist numbers $\lambda_i\in\mathbb{R}$ and a permutation $\pi\in S_n$ such that
    \[
        R(E_{ii}) = \lambda_i\otimes E_{\pi(i)\pi(i)}
        \quad\text{for all }1\leq i\leq n.
    \]
    Theorem~\ref{th_diag_k} gives the corresponding condition (i), (ii) or (iii) on $\lambda_i$ and~$\pi$. Set $\Lambda = \mathrm{diag}(\lambda_1,\ldots,\lambda_n)$, and let $P_2\in M_n$ be the matrix corresponding to the permutation $\pi^{-1}$. Then
    \[
        R(A) = N(A)\oplus (P_2\otimes\Lambda\otimes D(A)\otimes P_2^t) = T_5(A).
    \]
    
    The condition $\pi(i) = i$ is equivalent to the condition $(P_2)_{ii} = 0$. Hence the numbers~$\lambda_i$ and the matrix $P_2$ satisfy conditions (i)--(iii). The equality $T = T'\circ R$ gives the desired decomposition. The necessity is proved.

    3. Conversely, let the parameters from the statement of the theorem be given. Define the mapping $T'$ as $T_2 \circ T_3 \circ T_4$ or $T_1 \circ T_2 \circ T_3 \circ T_4$, according to the decomposition of $T$. Then $T = T'\circ T_5$.
    
    Let the matrix $P_2$ correspond to the permutation $\pi^{-1}$. It follows from the definition of the mapping $T_5$ that
    \[
    T_5(E_{ij}) = 
    \begin{cases}
        E_{ij}, & i \neq j,\\
        \lambda_i\otimes E_{\pi(i)\pi(i)}, & i = j.
    \end{cases}
    \]
    The mapping $T'$ is bijective and preserves the cyclicity index, as was shown above. The condition $(P_2)_{ii} = 0$ is equivalent to the condition $\pi(i) = i$. Then, by Theorem~\ref{th_diag_k}, the mapping $T_5$ preserves the value $k$ of the cyclicity index. Hence the composition $T = T'\circ T_5$ also preserves the value $k$ of the cyclicity index.
\end{proof}

\begin{corollary} \label{cor_preserving_k}
    Let $k \geq 3$, $n\geq\max(k, 4)$, and $T:M_n\to M_n$ be a bijective linear
    transformation of tropical matrices. If $n \geq m(k) + q(k)$ (in particular, by Proposition~\ref{proposition:case-(iii)}, this holds whenever $k$ is not a prime power and $k \neq 6$), then the following conditions are equivalent:
    \begin{enumerate}
        \item $T$ preserves the value $k$ of the cyclicity index;
        \item $T$ preserves the cyclicity index.
    \end{enumerate}
\end{corollary}

\begin{proof}
    The implication $(2)\Rightarrow (1)$ is obvious; we prove the implication $(1)\Rightarrow (2)$.
    
    Theorem~\ref{th_final_k} gives the decomposition $T = T' \circ T_5$, where $T'$ preserves the cyclicity index and $T_5$ satisfies the restriction (iii). Hence $T_5$ is the identity mapping, and $T = T'$ preserves the cyclicity index.
\end{proof}

\begin{corollary} \label{cor_preserving_k_both_sides}
    Let $k \geq 3$, $n\geq\max(k, 4)$, and $T:M_n\to M_n$ be a bijective linear transformation of tropical matrices. If $n > m(k)$ (in particular, this holds if $k = 3$ or, by Remark~\ref{remark:case-(i)}, whenever $k$ is not a prime power), then the following conditions are equivalent:
    \begin{enumerate}
        \item the mappings $T$ and $T^{-1}$ preserve the value $k$ of the cyclicity index,
        \item the mappings $T$ and $T^{-1}$ preserve the cyclicity index.
    \end{enumerate}
\end{corollary}

\begin{proof}
    The implication $(2)\Rightarrow (1)$ is trivial; we prove the implication $(1)\Rightarrow (2)$. If $n \geq m(k) + q(k)$, then the statement follows from Corollary~\ref{cor_preserving_k}. Assume that $m(k) < n < m(k) + q(k)$.
    
    Theorem~\ref{th_final_k} gives the decomposition $T = T' \circ T_5$, where $T'$ preserves the cyclicity index and $T_5$ satisfies the restrictions (i)--(ii). Thus $T_5^{-1} = T^{-1} \circ T'$. Since $T^{-1}$ preserves the value $k$ of the cyclicity index, $T_5^{-1}$ also does. One can verify that 
    \[
    T_5^{-1}(E_{ij}) = 
    \begin{cases}
        E_{ij}, & i \neq j,\\
        (-\lambda_{\pi^{-1}(i)})\otimes E_{\pi^{-1}(i)\pi^{-1}(i)}, & i = j.
    \end{cases}
    \]
    Theorem~\ref{th_diag_k} then implies that $-\lambda_{\pi^{-1}(i)} \leq 0$ for all $1 \leq i \leq n$, that is, all $\lambda_1, \dots, \lambda_n$ are zero. It follows from the restriction~(ii) that $P_2$ is the tropical identity matrix, and so $T_5$ is the identity mapping. Thus $T = T'$ and $T^{-1} = (T')^{-1}$ preserve the cyclicity index.
\end{proof}

\bibliographystyle{abbrvurl}
\bibliography{thebibliography-arXiv}

\end{document}